\documentclass{gtart}
\gtart

\usepackage{amsthm}
\usepackage{amsgen}
\usepackage{amsmath,amssymb}
\usepackage{curvesls}
\usepackage{epic}
\usepackage[all]{xy}
\usepackage[dvips]{graphicx}
\usepackage{psfrag}

\newtheorem{thm}{Theorem}[section]

\newtheorem{lem}[thm]{Lemma}

\newtheorem{prop}[thm]{Proposition}

\theoremstyle{definition}
\newtheorem{defn}[thm]{Definition}

\newcommand{\splice}{{\bowtie}}

\newcommand{\Real}{{\mathbb{R}}}

\newcommand{\Zed}{{\mathbb{Z}}}

\newcommand{\Diff}{{\mathrm{Diff}}}
\newcommand{\HomEq}{{\mathrm{HomEq}}}
\newcommand{\Isom}{{\mathrm{Isom}}}
\newcommand{\Trans}{{\mathrm{Trans}}}

\newcommand{\Emb}{{\mathrm{Emb}}}

\newcommand{\cubeact}{{\kappa}}

\newcommand{\EK}[1]{{\mathrm{EC}({#1})}}

\newcommand{\Cu}{{\mathcal{C}}}

\newcommand{\K}{{\mathcal{K}}}
\newcommand{\Prime}{{\mathcal{P}}}
\newcommand{\tK}{{\mathcal{\hat K}}}

\newcommand{\IG}{{\Bbb G}}

\begin{document}

\title{Topology of knot spaces in dimension $3$}

\author{Ryan Budney}

\address{
Mathematics and Statistics, University of Victoria\\
PO BOX 3045 STN CSC, Victoria, B.C., Canada V8W 3P4
}

\email{rybu@uvic.ca}

\begin{abstract} 
This paper gives a detailed description of the homotopy type of $\K$, the space of 
long knots in $\Real^3$, the same space of knots studied by Vassiliev via
singularity theory.  Each component of $\K$ corresponds to an isotopy class of
long knot, and we list the components via the companionship
trees associated to knots. The knots with the shortest companionship trees 
are: the unknot, torus knots, and hyperbolic knots.  
The homotopy type of these components of $\K$ were computed by Hatcher. 
In the case the companionship tree has more than one vertex, we give a fibre-bundle 
description of the corresponding components of $\K$, recursively, in terms of the 
homotopy types of components of $\K$ corresponding to knots with shorter
companionship trees. The primary case studied in this paper is the case of a knot
which has a hyperbolic manifold contained in the JSJ-decomposition
of its complement. Moreover, the homotopy type of $\K$ as an $SO_2$-space
is determined, which gives a detailed description of the homotopy-type of
the space of embeddings of $S^1$ in $S^3$.
\end{abstract}


\maketitle

\section{Introduction}\label{INTRODUCTION}

Let $\K$ denote the space of smooth embeddings
of $\Real$ in $\Real^3$ which agree with a fixed linear embedding outside of 
$[-1,1]$. $\K$ is called the space of long knots in $\Real^3$. The
fixed linear embedding will be called the `long axis.'  
To give a proper statement of the results in this paper, some familiarity
with the JSJ-decomposition of $3$-manifolds applied to knot complements
in $S^3$ as in the paper \cite{jsj} will be assumed. A very brief summary which
fixes some notation is included in the following four paragraphs.
 
Given a long knot $f \in \K$ the JSJ-decomposition of the knot complement 
$C$ consists of a collection
of incompressible tori $T \subset C$.  To be precise, consider $S^3$
to be the one-point compactification of $\Real^3$, thus $f$ has a one-point
compactification which is an embedding of $S^1$ in $S^3$.  $C$ is the
complement of an open tubular neighbourhood of this embedded
$S^1$ in $S^3$. Let $C | T$ denote the complement of an open tubular
neighbourhood of $T$ in $C$.
Form a graph $G_f$ whose vertices are the path-components of $C | T$, and whose 
edges are path-components of $T$.  This graph has the structure of a 
rooted tree, with root $V$ given by the component of $C | T$ containing 
$\partial C$.  Call this tree the JSJ-tree of $f$. The root manifold
$V$ has the interesting property that it has a canonical re-embedding
in $S^3$ as the complement of a Knot Generating Link (KGL) -- this is a link
$L=(L_0,L_1,\cdots,L_n)$ whose complement is prime and either
Seifert fibred or atoroidal, such that the $n$-component
sub-link $(L_1,L_2,\cdots,L_n)$ is the unlink.  $n$ is the number of
children of $V$ in the JSJ-tree. Moreover, 
$\overline{C \setminus V} = C_{J_1} \sqcup \cdots \sqcup C_{J_n}$
is the disjoint union of $n$ non-trivial knot complements.
Thus, the knot $f$ is can be represented as a `spliced' knot, $f=J \splice L$
where $J=(J_1,\cdots,J_n)$ is an $n$-tuple
of non-trivial long knots whose complements are $C_{J_1}, \cdots, C_{J_n}$. One 
obtains the spliced knot $J \splice L$ from $L_0$ by a
generalized satellite construction: fix
$D = (D_1,\cdots,D_n)$ a collection of disjoint discs in $S^3$ 
such that $\partial D = (\partial D_1,\cdots,\partial D_n) = (L_1,\cdots,L_n)$,
and fix $\nu D$ a closed tubular neighbourhood of $D$.
$J \splice L$ is the union of $L_0 \setminus \nu D$  together with the (untwisted) 
re-embedding of $L_0 \cap \nu D_i$ along the knots $J_i$ for all 
$i \in \{1,2,\cdots,n\}$. Specifically, re-embed a tubular neighbourhood of 
$D_i$ using $J_i$ as the `pattern,' like in the diagram below.

{
\psfrag{f}[tl][tl][1][0]{$f$}
\psfrag{L}[tl][tl][1][0]{$L$}
\psfrag{J1}[tl][tl][1][0]{$J_1$}
\psfrag{J2}[tl][tl][1][0]{$J_2$}
$$\includegraphics[width=9cm]{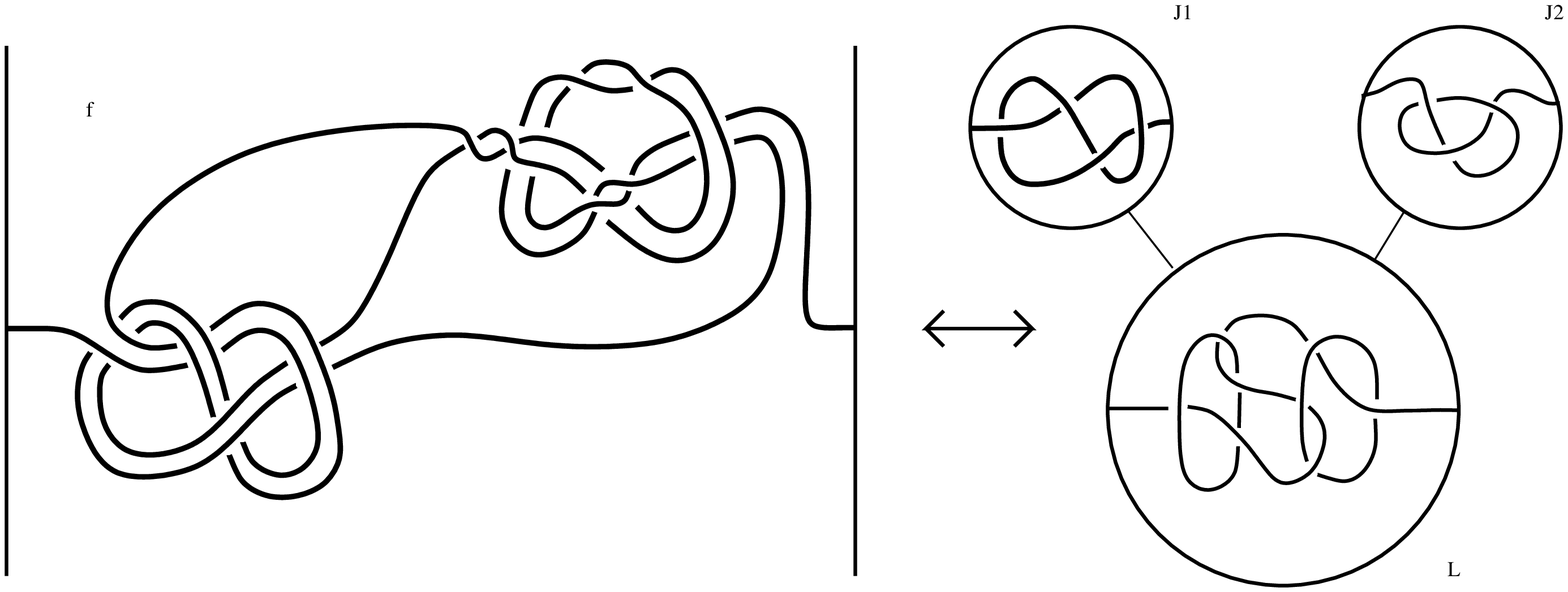}$$
\centerline{$f = J \splice L$ where $L$ is the Borromean rings}
}

When $n=0$, interpret $J \splice L$ to be $L$. In this case $L$ is either
a torus knot, a hyperbolic knot or an unknot. In the general case, after one 
writes $f$ as a splice, $f=J \splice L$, one can
continue and write the knots of $J$ as spliced knots, iteratively.
This process terminates, at which point the iterated splice construction
for $f$ is precisely the JSJ-tree of $C$, with every vertex
decorated by the appropriate KGL. This decorated
tree $\IG_f$ is called the companionship tree of $f$. $\IG_f$ 
is a complete isotopy invariant of long knots.  For the purpose
of this paper, companionship trees will be used to index the components of $\K$.

Seifert-fibred KGLs were classified by Burde and Murasugi. 
Thurston proved that non-Seifert fibred KGLs are hyperbolic.
Kanenobu proved that hyperbolic KGLs exist with arbitrarily many
components, thus any abstract finite rooted tree can be
realized as the JSJ-tree $G_f$ of some long knot $f$. 
See \cite{jsj} for complete references to all the above comments
regarding JSJ-decompositions and companionship trees. 

Given a knot $f \in \K$, label the path-component of $\K$ containing $f$
by $\K_f$. As a set, this is all the long knots isotopic to $f$.  
The primary result of this paper is the computation of the homotopy type of
$\K_f$ if $f$ is a hyperbolically-spliced knot ie: $f=J \splice L$ where
$L$ is a hyperbolic KGL. 
Let $\Sigma_n^+$ be the signed symmetric group. $\Sigma_n^+$ is the 
semi-direct product $\Zed_2^n \rtimes \Sigma_n$ where $\Sigma_n$ acts
on $\Zed_2^n$ by permutation of factors. 
Given $L = (L_0,L_1,\cdots,L_n)$ a hyperbolic KGL, let $B_L$ be the group of 
hyperbolic isometries of the complement of
$L$ which extend to orientation-preserving diffeomorphisms of $S^3$ 
such that the extensions to $S^3$ preserve $L_0$ and its orientation.
$B_L$ is known to be cyclic by the proof of the Smith Conjecture \cite{SmithConj}.
There is a natural representation $B_L \to \Sigma^+_n$ given
by the composite 
$B_L \to \Diff(S^3,L) \to \pi_0 \Diff(\cup_{i=1}^n L_i) \equiv \Sigma_n^+$.  
There is an action of $\Sigma_n$ on $\K^n$ given by
permutation of the factors. There is an action of $\Zed_2$ on $\K$ given
by knot inversion, meaning, rotation by $\pi$ about a fixed axis perpendicular 
to the long axis.  These two actions extend in a unique way 
to an action of $\Sigma_n^+$ on $\K^n$. Thus via the representation 
$B_L \to \Sigma_n^+$ there is an induced action of $B_L$ on $\K^n$.
Let $A_f$ be the maximal subgroup of $B_L$ such that the induced action of 
$A_f$ on $\K^n$ preserves $\prod_{i=1}^n \K_{J_i}$.

\begin{thm}{\bf (Proved as Theorem \ref{mainthm}) }
\label{previewthm}
If $f=J \splice L$ where $L$ is an $(n+1)$-component hyperbolic KGL, then 
$$\K_f \simeq S^1 \times \left( SO_2 \times_{A_f} \prod_{i=1}^n \K_{J_i} \right)$$
The composite $A_f \to \Diff(S^3,L_0) \to \Diff(L_0)$
is faithful, giving an embedding $A_f \to SO_2$, which is the action
of $A_f$ on $SO_2$.
\end{thm}

This result gives a recursive description of all of the homotopy-type of
all of the components of $\K$ since we have the prior results:
 
\begin{enumerate}
 \item \cite{Hatcher2} If $f$ is the unknot, then $\K_f$ is contractible.
 \item \cite{Hatcher4} If $f$ is a torus knot, then $\K_f \simeq S^1$.
 \item \cite{Hatcher4} If $f$ is a hyperbolic knot, then 
  $\K_f \simeq S^1 \times S^1$.
 \item \cite{Hatcher4} If a knot $f$ is a cabling of a 
 knot $g$ then $\K_f \simeq S^1 \times \K_g$. 
 \item \cite{cubes} If the knot $f$ is a connected sum of
 $n \geq 2$ prime knots $f_1, f_2, \cdots, f_n$ then
 $\K_f \simeq \Cu_2(n) \times_{\Sigma_f} \prod_{i=1}^n \K_{f_i}$.
 $\Cu_2(n) = \{ (x_1,\cdots,x_n) \in (\Real^2)^n : x_i \neq x_j \ \forall i \neq j \}$.
 $\Sigma_f \subset \Sigma_n$ is a Young subgroup of $\Sigma_n$, acting on
 $\Cu_2(n)$ by permutation of the coordinates, and similarly by permuting
 the factors of the product $\prod_{i=1}^n \K_{f_i}$.  $\Sigma_f$ is the subgroup
 of $\Sigma_n$ preserving the partition of $\{1,2,\cdots,n\}$ defined by the equivalence
 relation $i \sim j \Longleftrightarrow \K_{f_i} = \K_{f_j}$.
 \item \cite{jsj} If a knot has a non-trivial companionship tree, then
  it is either a cable, in which case (4) applies, 
  a connect-sum, in which case (5) applies, or is hyperbolically
  spliced, in which case Theorem \ref{previewthm} applies. If
  a knot has a trivial companionship tree, it is either the unknot,
  in which case (1) applies, or a torus knot in which case
  (2) applies, or a hyperbolic knot, in which case (3) applies.
  Moreover, every time one applies one of the above theorems,
  one reduces the problem of computing the homotopy type of 
  $\K_f$ to computing the homotopy type of knotspaces for knots 
  with shorter companionship trees, thus the process terminates.
\end{enumerate}

Items (2) and (4) are re-proved in Theorem \ref{cable_thm},
(3) is a special case of Theorem \ref{previewthm}, and 
(5) is re-proven in Theorem \ref{sumthm}.

In Section \ref{htpy_r3} Theorem \ref{mainthm} is proved, giving the homotopy
type of $\K_f$ when $f = J \splice L$, $L$ a hyperbolic KGL, completing
the description of the homotopy type of $\K$, modulo one issue.
That being, a recursive description of the homotopy class of the
inversion map $I : \K_f \to \K_f$ when $f$ is an invertible knot. 

In Section \ref{3sphere} a summary is given of how the results of 
Sections \ref{cables} and \ref{htpy_r3} enable one to compute the homotopy 
type of $\K$ as an $SO_2$-space, where the $SO_2$ is the subspace of 
$SO_3$ given by rotations that preserve the `long axis.'
In Theorem \ref{eq3} a recursive computation of the homotopy class of the inversion map
is given.  These results in turn describe the homotopy type of
the space of embeddings of $S^1$ in $S^3$.

In Section \ref{Examples} various algorithms and computational
methods are mentioned that allow one to apply the theorems of
this paper and compute the homotopy type of components of $\K$.
Several example computations using SnapPea are given.

Section \ref{gramelt} is devoted to `the Gramain element.' This is
a non-trivial element of $\pi_1 \K_f$ provided $f$ is non-trivial. 
It was first studied by Gramain \cite{Gramain2}, and is
given by $2\pi$ rotation around the long axis. The Gramain element 
survives to give a non-trivial element in $H_1 \K_f$ which is
called the Gramain cycle. Thus, the only component of $\K$ with trivial 
homology is the unknot component. 

Section \ref{endgame} considers several natural questions that the
results of this paper raise, such as: what is the relationship to the 
Vassiliev spectral sequence for $H_* \K$, the Goodwillie calculus,
the homotopy type of the group-completion of $\K$, 
and a symmetry realization problem for hyperbolic KGLs related
to the problem of giving a `closed form' description of $H_* (\K;\Zed)$.

%%% namedropping %%%

The study of embeddings and to a lesser extent, spaces of embeddings, was initiated 
by Whitney \cite{Whitney}, whose
Weak Embedding Theorem could perhaps be viewed as one of the first theorems in
the subject.  Of course, classical knot theory is an example of the theory
of embedding spaces and it predates Whitney, but the focus of the subject historically
was not on the global nature of embedding spaces, but 
more on the `discrete' isotopy classes of knots. 
 
Perhaps the first `strong' results on the homotopy type of $\K$ were due 
to Gramain \cite{Gramain2}, who used the Burde-Zieschang
characterization of the centres of knot groups to show that the 
map $\pi_1 SO_2 \to \pi_1 \K_f$ given by rotation of a non-trivial long 
knot about the `long axis' is injective.  Gramain also constructed
injections $P_n \to \pi_1 \K_f$ for certain non-prime $f \in \K$
where $P_n$ is the pure braid group.

By the work of Hatcher \cite{Hatcher1} and Ivanov \cite{Ivanov1, Ivanov2}
it was known that the Smale Conjecture is equivalent to 
proving that the unknot component of $\K$ is contractible, 
which Hatcher later went on to prove \cite{Hatcher2}.
Less well known is one of the consequences: that all components of
$\K$ are $K(\pi,1)$'s \cite{Hatcher4}.  
More recently, Hatcher and McCullough proved \cite{HM} that 
each component of $\K$ has the homotopy type 
of a finite-dimensional CW-complex. 
In a sequence of preprints \cite{Hatcher4, Hatcher5}
Hatcher has also shown, assuming the Linearization Conjecture
(that every finite group acting on $S^3$ is conjugate to a linear
action), how to compute the homotopy type of $\K_f$ provided the companionship
tree $\IG_f$ has no branching.  Moreover, Hatcher went on
to describe the homotopy type of the components of hyperbolic knots
and torus knots in $\Emb(S^1,S^3)/\Diff(S^1)$.

For the author, it was Victor 
Tourtchine's conjecture that knot spaces admit actions of the operad
of little $2$-cubes that motivated recent work in the subject.  After showing that
$\K$ admits a little $2$-cubes action \cite{cubes} it was only natural to ask
if the action was `free' in the sense of May \cite{May}. The affirmative answer to
that question \cite{cubes} along with some work on JSJ-decompositions \cite{jsj} led to
the realization that the homotopy type of $\K$ is quite explicitly `computable,'
which is the point of this paper.

Tourtchine \cite{Turchin2} was working on the homotopy type of spaces of knots 
from the point of view of singularity theory, as developed by Vassiliev \cite{Vass}.
He described a bracket on the $E^2$-page of the Vassiliev spectral sequence for $H_* \K$ which conjecturally came from an action of the operad of 
little $2$-cubes on $\K$.  Vassiliev's theory looks at $\K$ as a subspace of a 
contractible space of maps, $\K$ being seen as the complement of a 
`discriminant' space of singular maps, this gives rise to
 a spectral sequence for the homology of $\K$ via a filtration of the 
 `discriminant' by a natural notion of `complexity' of the singularity.
 Vassiliev's spectral sequence has been, in effect, re-derived by Sinha \cite{Dev}
 who studies spaces of knots using the Goodwillie Calculus of Embeddings \cite{Good}. 
Both the approaches of singularity theory and Goodwillie Calculus have
technical problems in dimension $3$ as the relevant spectral sequences have 
been shown to not converge to $H_* \K$ \cite{Vass, Volic}. 
Indeed, most of the literature related to Vassiliev's
work focuses on the elements of $H^0 \K$ produced by the cohomological 
spectral sequence, the so-called `Vassiliev invariants' of knots. What relation 
the results in this paper have to these spectral-sequence approaches to $\K$ is an 
open problem.

The Goodwillie Calculus methods work very well in
co-dimension 3 or larger, as they allow explicit computation of various
homotopy groups of knotspaces \cite{DevKevin}. Recently Lambrechts,
Tourtchine and Volic \cite{Lambrechts} have reduced the problem
of computing the rational homology of such knot
spaces to computing the homology of a rather explicit DGA.
They do this by showing the rational collapse of the 
Vassiliev spectral sequence. There are several new interesting 
developments in the area, such as K.~Sakai's homological 
relationship between $\Omega C_k(\Real^n)$ and the space of long 
knots in $\Real^{n+1}$ for $n \geq 3$ via `capping.'  In a rather 
direct application of Goodwillie's dissertation, it is shown that 
the first non-trivial homotopy groups of the space of `long 
embeddings' of $\Real^j$ in $\Real^n$ are related via a generalized 
Litherland Spinning construction provided $2n-3j-3 \geq 0$ \cite{embfam}. 
Among other things this result provides a geometric isomorphism between 
$\pi_2(\text{the space of
long knots in } \Real^4)$ and $\pi_0 \Emb(S^3,S^6) \simeq \Zed$.

%%% /namedropping %%%
 
I thank Allen Hatcher for his input on this
project. I thank Fred Cohen.  This paper probably never 
would have been finished if it wasn't for his frequent
encouragement. I'd like to thank Cameron Gordon for pointing out an error in 
an early attempt at a proof of Theorem \ref{mainthm} and for bringing the
Hartley paper to my attention \cite{Hartley}. 
I owe thanks to Jeff Weeks, indirectly, for his wonderful software, SnapPea.  
I thank Morwen Thistlethwaite for showing me how to convince 
SnapPea to compute the symmetry groups $B_L$ 
in Theorem \ref{mainthm}. I thank Yoshiaki Uchida for his input.  
I also thank the University of Oregon's topology/geometry seminar, 
particularly Dev Sinha, Hal Sadofsky, and Gregor Masbaum for listening
to early versions of these results, and the resulting helpful conversations.
I would also like to thank the Max Planck Institute for Mathematics and
Institut des Hautes \'Etudes Scientifiques for support during the preparation
of this manuscript. 

\section{The recursive setup for the homotopy-type of $\K$}\label{cables}

This section reviews the techniques for describing the homotopy-type of 
$\K$, quoting and reproving various theorems about embedding spaces. 

\begin{thm}\label{fundprop}
Given a long knot $f \in \K$, consider the one-point compactification of
$f$ to be an embedding of $S^1$ in $S^3$ and let $C$ be the complement
of an open tubular neighbourhood of this embedding. The restriction
map is a Serre fibration 
$$ \Diff(D^3 \text{ fix } \partial D^3) \to \K_f $$
where $\Diff(D^3 \text{ fix } \partial D^3)$ is the group of diffeomorphisms of 
$D^3$ which fix the boundary disc pointwise.  
\begin{enumerate}
%\item $\Diff(D^3) \to \K_f$ is null-homotopic.
\item \cite{Hatcher4} The fibre of the map $\Diff(D^3 \text{ fix } \partial D^3) \to \K_f$ 
has the homotopy type of $\Diff(C \text{ fix } \partial C)$, the group of diffeomorphisms of 
$C$ that fix the boundary pointwise. 
%\item There is an injection $\pi_{i+1} \K_f \to \pi_i \Diff(C)$
%for all $i$ and the cokernel is isomorphic to $\pi_i \Diff(D^3)$.
\item \cite{Hatcher2} $\Diff(D^3 \text{ fix } \partial D^3)$ is contractible, thus 
$\K_f \simeq B\Diff(C \text{ fix } \partial C)$.
\item \cite{Hatcher1} $\Diff(C \text{ fix } \partial C)$ has the homotopy-type of a discrete 
group. Equivalently, $\K_f$ is a $K(\pi,1)$.
\end{enumerate}
\end{thm}

Given any manifold $M$ and a submanifold $N$, that the restriction map 
$\Diff(M) \to \Emb(N,M)$ is a Serre fibration is a classical theorem.  The homotopy 
lifting property for this map is proved by observing that in the classical proof of the 
Isotopy Extension Theorem for an isotopy of $N$ in $M$, the solution depends smoothly on 
the initial conditions (the embedding of $N$ in $M$), as the isotopy extension is
constructed by solving an ordinary differential equation. 
A specific reference is Cerf's paper \cite{Cerf} \S 2.2.2, but perhaps it is preferable
for the reader to adapt Hirsch's textbook proof of isotopy extension \cite{Hirsch} \S 8.1.3.
Point (1) in Proposition \ref{fundprop} follows from the homotopy classification
of spaces of tubular neighbourhoods, see \cite{Cerf} or adapt the 
isotopy-classification of tubular neighbourhoods from \cite{Hirsch}. Point (2)
is a difficult theorem of Hatcher's, historically called `the Smale Conjecture' 
\cite{Hatcher2}, and point (3) is a consequence of (2) when combined with the
theorems, proved independently by Hatcher and Ivanov, on the homotopy-type of
spaces of incompressible surfaces in 3-manifolds \cite{Hatcher1, Ivanov1, Ivanov2}.

\begin{lem}\label{mainlemma}Let $f$ be a compound knot, ie: 
$f = (f_1,\cdots,f_n) \splice L$ with $L$ a KGL. 
There is an exact sequence of groups
$$0 \to \prod_{i=1}^n \pi_0 \Diff(C_{f_i} \text{ fix } \partial C_{f_i}) \to 
    \pi_0 \Diff(C \text{ fix } \partial C) \to \pi_0 \Diff(V \text{ fix } \partial C)$$
where $\overline{C \setminus V} = \sqcup_{i=1}^n C_{f_i}$, ie: $C_{f_i}$ is diffeomorphic
to the complement of a non-trivial knot $f_i$.  $V$ is canonically diffeomorphic to the complement
of the link $L$, denoted $C_L$.  
\begin{proof}
The key technique is to observe that $\Diff(C \text{ fix } \partial C)$ deformation-retracts
to the subgroup which preserves $V$, using the incompressible surface results in \cite{Hatcher1}
or \cite{Ivanov1, Ivanov2}.  Restriction to $V$ gives a fibration whose homotopy long
exact sequence is the exact sequence in the lemma. 
\end{proof}
\end{lem}

Since there is a canonical diffeomorphism $V \to C_L$ (called the untwisted re-embedding in \cite{jsj}), 
let $T_0 \subset \partial C_L$ correspond to $\partial C \subset V$ under this diffeomorphism. 
Observe that image of the map 
$\pi_0 \Diff(C \text{ fix } \partial C) \to 
 \pi_0 \Diff(V \text{ fix } \partial C) \equiv \pi_0 \Diff(C_L \text{ fix } T_0)$ 
consists precisely of the elements of $\Diff(C_L \text{ fix } T_0)$ which extend
to diffeomorphisms of the pair $(S^3,L)$ such that
(a) if $L_i$ is sent to $L_j$ preserving the orientations, then 
$f_i$ is isotopic to $f_j$, and
(b) if $L_i$ is sent to $-L_j$, then $f_i$ is isotopic to the 
inverse of $f_j$.  The conventions of Lemma \ref{mainlemma} will be used throughout
the remainder of the paper.

\begin{thm}\label{cable_thm}
If $f\in \K$ is a cable of $g \in \K$ then $\K_f \simeq S^1 \times \K_g$.
\begin{proof}
There are two cases to consider:
\begin{enumerate}
\item The JSJ-tree has only one vertex.  This is the case where
$f$ is a non-trivial torus knot, $g$ is the unknot, and $V=C$. 
\item The JSJ-tree has two or more vertices.  This is the case that
$f$ is a cable of a non-trivial knot $g$ and $V$ is diffeomorphic to
the complement of a $2$-component Seifert-fibred link complement in
$S^3$. 
\end{enumerate}
In case (1), $\pi_0 \Diff(C \text{ fix } \partial C) \simeq \Zed$.  This is because
$C$ is a Seifert-fibred manifold, fibred over a disc with two singular
fibres.  The singular fibres have Seifert-fibre data $\frac{\alpha_1}{\beta_1}$ and 
$\frac{\alpha_2}{\beta_2}$ where $\alpha_1\beta_2 + \alpha_2\beta_1=\pm 1$, 
in the `unnormalized' notation of Hatcher $C \simeq M(0,1;\frac{\alpha_1}{\beta_1},
\frac{\alpha_2}{\beta_2})$ \cite{Hatcher3}.
Any diffeomorphism of $C$ can be isotoped to one which preserves
the Seifert-fibring \cite{Hatcher1}, thus $\pi_0 \Diff(C)$ is generated by
a meridional Dehn twist about a boundary-parallel torus.  Such Dehn twists are
non-trivial since all non-trivial powers of the meridian are non-central in $\pi_1 C$
\cite{BZ}, thus $\K_f \simeq B\Diff(C) \simeq S^1$ and $\pi_1 \K_f$ is generated
by full rotation about the long axis of the knot.
{
\psfrag{pq}[tl][tl][1.2][0]{$\frac{\alpha_1}{\beta_1}$}
\psfrag{qp}[tl][tl][1.2][0]{$\frac{\alpha_2}{\beta_2}$}
$$\includegraphics[width=7cm]{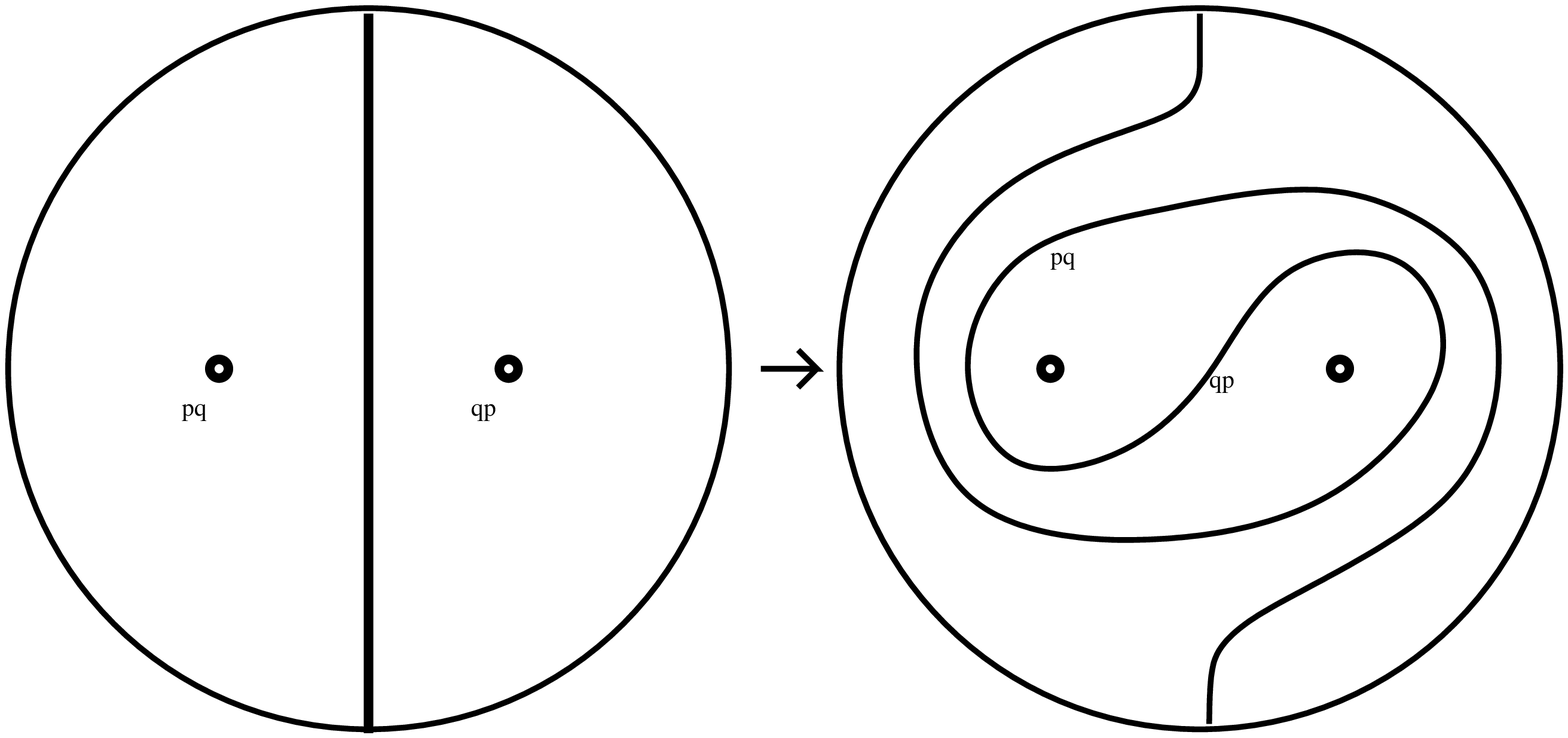}$$
\centerline{The generator of $\pi_0 \Diff(C) \simeq \pi_1 \K_f$}
}
In case (2) $V \simeq M(0,2;\frac{\alpha_1}{\beta_1})$ 
is Seifert fibred over a once-punctured disc with one singular fibre. 
The map $\pi_0 \Diff(C \text{ fix } \partial C) \to \pi_0 \Diff(V \text{ fix } \partial C)$
is onto an infinite-cyclic group generated by a meridional Dehn twist about
a torus parallel to $\partial C$ via an argument analogous to case (1). This allows the 
direct construction of an isomorphism
$\pi_0 \Diff(C_g \text{ fix } \partial C_g) \times \Zed \to \pi_0 \Diff(C \text{ fix } \partial C)$
by taking the union of a diffeomorphism on $C_g$ with a Dehn twist on $V$.
\end{proof}
\end{thm}

\begin{thm}\label{sumthm}If $f$ is a connect sum of prime knots $f_1, \cdots, f_n$ then
$$\K_f \simeq \Cu_2(n) \times_{\Sigma_f} \prod_{i=1}^n \K_{f_i}$$
where $\Sigma_f \subset \Sigma_n$ is the Young subgroup corresponding
to the equivalence relation $i \sim j$ iff $f_i$ and $f_j$ are isotopic
as long knots.
\begin{proof}
The proof follows the outline of Theorem \ref{cable_thm}, but in this case
$V \simeq C_L$ where $L$ is the $(n+1)$-component keychain link.  
The mapping class group of diffeomorphisms of the pair $(S^3,L)$ which fix a tubular 
neighbourhood of $L_0$ pointwise is isomorphic to the braid group, since
$C_L$ is diffeomorphic to a product of a circle with an $n$-punctured disc \cite{Hatcher1}.
Thus the image of the map $\pi_0 \Diff(C \text{ fix } \partial C) \to \pi_0 \Diff(C_L \text{ fix } T_0)$
corresponds to the preimage of $\Sigma_f$ under the canonical epi-morphism $B_n \to \Sigma_n$
where $B_n$ is the braid group on $n$-strands.  Constructing a splitting over this
group amounts to choosing a canonical identification between
$C_{f_i}$ and $C_{f_j}$ for all $f_i$ isotopic to $f_j$ and extending via this
identification.
\end{proof}
\end{thm}

In \cite{cubes} it took quite a bit more work to prove Theorem \ref{sumthm}
since it involved showing that an action of the operad of $2$-cubes gave the splitting.

\section{The homotopy type of a hyperbolically spliced knot}\label{htpy_r3}

Let $V$ be an orientable, complete cusped hyperbolic manifold of finite volume.  
Let $\overline{V}$ denote the manifold compactification of $V$, ie: if $V$ has
$n$ cusps, then $\overline{V}$ has $n$ torus boundary components. 

\begin{defn}
Let $\Isom(V)$ be the group of hyperbolic isometries of $V$. 
Any isometry of $V$ extends to a diffeomorphism of $\overline{V}$ since
in a horoball neighbourhood of a cusp the isometry is a Euclidean motion.
This gives a map $\Isom(V) \to \Diff(\overline{V})$.
Let $\HomEq(V)$ be the space of homotopy equivalences of $V$. There is 
a natural map $\Diff(\overline{V}) \to \HomEq(V)$ given by restriction to $V$.
\end{defn}

\begin{prop}\label{hyp-smooth}
The maps $\Isom(V) \to \Diff(\overline{V})$ and $\Diff(\overline{V}) \to \HomEq(V)$
are homotopy equivalences.
\begin{proof}
Mostow rigidity  \cite{Mostow} tells us that the map $\Isom(V) \to \HomEq(V)$ is a homotopy 
equivalence. By the work of Hatcher and Ivanov \cite{Hatcher1, Ivanov1, Ivanov2}, $\Diff(\overline{V})$ is a 
$\K(\pi,0)$. Thus we need only show that $\pi_0 \Diff(\overline{V}) \to \pi_0 \HomEq(V)$ is 
injective.  Waldhausen \cite{Wald} (Lemma 1.4.2) proves homotopic diffeomorphisms of 
Haken manifolds are isotopic provided the homotopy 
preserves the boundary of the manifolds. In his Theorem 7.1
\cite{Wald} he proves this boundary condition can always 
be satisfied for boundary-irreducible manifolds.  We can disregard the special case
of $I$-bundles as they are not cusped hyperbolic manifolds of finite volume.
\end{proof}
\end{prop}

Let $L=(L_0,L_1,\cdots,L_n)$ be an $(n+1)$-component hyperbolic KGL, and let
$C_L$ be its complement. Let $T_0$ be the boundary of $C_L$ corresponding to the $L_0$-cusp. 
Note that an isometry of $C_L$ extends to a diffeomorphism of $S^3$ if and only if it
preserves the longitudes and meridians of $L$.  Moreover, if $G$ is any subgroup of
$Isom(C_L)$ consisting of elements which extend to diffeomorphisms of $S^3$, then the action
of $G$ on $C_L$ extends to an action of $G$ on the pair $(S^3,L)$. This action can be
constructed `by hand' using the fact that $G$ acts by Euclidean motions on the boundary
tori of $C_L$.

\begin{defn}
Let $B_L$ be the subgroup of $\Isom(C_L)$ that both extend to diffeomorphisms of $S^3$ and act on 
the $L_0$-cusp by translations, where we identify $\partial C_L$ with the boundary of some tubular
neighbourhood of $L$ in $S^3$.  Let $\Trans(T_0)$ be the group of Euclidean translations of
the $L_0$-cusp $T_0$. Basic Riemmanian geometry tells us that $B_L \to \Trans(T_0)$ is an embedding. 
As a group $\Trans(T_0)$ is isomorphic to $SO_2^2$. Let $\widetilde{\Trans}(T_0) \to \Trans(T_0)$
be the universal covering space, and let $\widetilde{B_L}$ be the lift of $B_L$ to the universal cover.
Let $\Diff_{T_0}(C_L)$ be the group of diffeomorphisms of $C_L$ which both extend to diffeomorphisms of 
$S^3$ and restrict to diffeomorphisms of $T_0$ which are isotopic to the identity. 
Let $\Diff(C_L \text{ fix } T_0)$ be the subgroup of $\Diff_{T_0}(C_L)$ which 
fixes $T_0$ pointwise.
\end{defn}

\begin{prop}\label{fractional-twists}
$\pi_0 \Diff(C_L \text{ fix } T_0)$ is free abelian on two generators.
Moreover the map $\Diff(C_L \text{ fix } T_0) \to \Diff_{T_0}(C_L)$ is equivalent to
the projection $\widetilde{B_L} \to B_L$.
\end{prop}

Proposition \ref{fractional-twists} is proven as Proposition 3.2 in \cite{HM}, so
the proof is omitted here.  But it is enlightening to understand the equivalence
$\xi : \widetilde{B_L} \to \Diff(C_L \text{ fix } T_0)$.  If $\tilde f \in \widetilde{B_L}$, 
let $f \in \Isom(C_L)$ be the corresponding isometry, then 
$\xi(\tilde f) \in \Diff(C_L \text{ fix } T_0)$ is defined to
be equal to $f$ outside a fixed collar neighbourhood of $T_0$.  Extend $\xi(\tilde f)$ to 
be a diffeomorphism of $C_L$ by the straight line from $0 \in \widetilde{\Trans}(T_0)$ to 
$\tilde f \in \widetilde{\Trans}(T_0)$. Such a path gives a path of translations of $T_0$ thus a 
fibre-preserving diffeomorphism of the collar neighbourhood of $T_0$ starting at the identity and
ending at a map which agrees with the isometry $f$.

For $L = (L_0,L_1,\cdots,L_n)$ a KGL, let $\Diff(S^3,L,L_0)$ be the group of diffeomorphisms
of $S^3$ which preserve $L$, and which also preserve $L_0$ together with its orientation.

\begin{prop}\label{badmono}{\bf (No Bad Monodromy)}
Let $L$ be a hyperbolic KGL, then  
the composite $B_L \to \Diff(S^3,L,L_0) \to \Diff(L_0)$ has a trivial kernel.
Equivalently, $B_L \to \Diff(L_0)$ is an embedding.
\begin{proof}
Let $H$ be the kernel of the composite $B_L \to \Diff(S^3,L,L_0) \to \Diff(L_0)$.
For $n=0$ this is true since if $H$ was non-trivial, $L = (L_0)$ is a hyperbolic 
knot that's a fixed point set of a finite group action on $S^3$ which would contradict the
Smith conjecture \cite{SmithConj}.
 
For $n \geq 1$ if $H$ was non-trivial, $L_0$ would be the unknot and $S^3/H$ can be identified 
with $S^3$ by the resolution of the Smith conjecture \cite{SmithConj}.  The quotient 
$(S^3,L)/H$ is a manifold pair $(S^3/H,L/H)$ where $S^3/H \simeq S^3$ and 
$L/H \subset S^3/H$ is some link. We identify $S^3/H$ with $S^3$ and $L/H$ with
$L'=(L_0',L_1',\cdots,L_p') \subset S^3$. Notice that $S^3 \setminus L'$ is hyperbolic.

By definition, there is an epimorphism 
$\pi_1 (S^3 \setminus L_0') \to H$.
Define a homomorphism 
$\pi_1 (S^3 \setminus L') \to \Zed_2$ to be the sum of the mod-$2$
linking numbers of a path in $S^3 \setminus L'$ with the link
components $L'\setminus L'_0$. Thus we have an epimorphism
$\pi_1 (S^3 \setminus L') \to \Zed_2 \oplus H$ given by the 
direct-sum of the above two homomorphisms.  

If we take the abelian
branched cover $\pi : X \to S^3$ corresponding to the above 
homomorphism $\pi_1 (S^3 \setminus L') \to \Zed_2 \oplus H$,
the branch-point set is precisely $L'$, and since the
homomorphism $\pi_1 (S^3 \setminus L') \to \Zed_2 \oplus H$
is a direct sum, the covering map 
$(X,\pi^{-1} L') \to (S^3,L')$ factors as a composite
$(X,\pi^{-1} L') \to (S^3,L) \to (S^3,L')$
where the map $(X,\pi^{-1} L') \to (S^3,L)$ is an abelian
branched cover with branch-point set $L \setminus L_0$.

Since $L \setminus L_0$ is an $n$-component unlink, 
this forces $X$ to be a connect-sum of copies of
$S^1 \times S^2$.  Precisely, 
$$X \simeq \#_{n-1} S^1 \times S^2$$
Sakuma's Theorem B \cite{Sakuma} tells us that $L'$ is either a trivial knot or
a connect-sum of Hopf links, or is split. Since $L'$ is hyperbolic it can
be none of these, contradicting our assumption that $H$ was non-trivial.
\end{proof} 
\end{prop}

\begin{defn}\label{def1}
Let $I : \K \to \K$ be the involution of $\K$ induced by rotation by 
$\pi$ around an axis perpendicular to the long axis. 
A knot $f \in \K$ is said to be invertible if $I(\K_f)=\K_f$. 
The knot $f$ is strongly invertible if $I : \K_f \to \K_f$ has a fixed point.
\end{defn}

Not all knots are invertible \cite{Trotter}, and not all invertible knots are
strongly invertible \cite{Hartley, Boileau}.
$\Zed_2 = \{\pm 1\}$ acts on $\K$ by inversion (rotation by
$\pi$ about any axis perpendicular to the long axis).  $\Sigma_n$ acts on
$\K^n$ by permutation of the factors.  These two actions induce an
action of $\Sigma_n^{+} = \{\pm 1\}^n \rtimes \Sigma_n$ on $\K^n$.
Let $f=J \splice L$ where $L$ is an $(n+1)$-component hyperbolic KGL, then
$B_L$ has a natural representation into $\Sigma_n^+$ given by
the composite $B_L \to \Diff(S^3,L,L_0) \to \pi_0 \Diff(L_1 \cup \cdots \cup L_n) \equiv \Sigma_n^+$.

\begin{defn}
$A_f$ is the subgroup of $B_L$ that preserves $\prod_{i=1}^n \K_{J_i}$ under the $B_L$ action
on $\K^n$.
\end{defn}

\begin{thm}\label{mainthm}
Let $f = J \splice L$, where $L$ is a hyperbolic KGL, 
then there is a homotopy-equivalence:
$$\K_f \simeq \left(\prod_{i=1}^n \K_{J_i}\right) \times_{A_f} \Trans(T'_0)$$
and there is a further splitting
$$\K_f \simeq S^1 \times \left( \left(\prod_{i=1}^n \K_{J_i}\right)
 \times_{A_f} L_0 \right)$$
\begin{proof} 
Let $\tilde A_f \to A_f$ be the pull-back of $\tilde B_L \to B_L$
via the inclusion $A_f \subset B_L$.  By Proposition \ref{hyp-smooth} and
the observations in Lemma \ref{mainlemma} 
the map $\pi_0 \Diff(C \text{ fix } \partial C) \to \pi_0 \Diff(C_L \text{ fix } T_0) \equiv \tilde B_L$
is onto $\tilde A_f$. By Proposition \ref{badmono}, $\tilde A_f$ is free abelian, generated by two elements, 
with one being a meridional Dehn Twist.  This allows us to construct a splitting of the
map $\pi_0 \Diff(C \text{ fix } \partial C) \to \tilde A_f$ in two steps: 1) choose an element
in the preimage of the generator which is not the meridional Dehn Twist.  2) The meridional Dehn Twist in 
$\pi_0 \Diff(C_L \text{ fix } T_0)$ has a meridional Dehn twist in 
$\pi_0 \Diff(C \text{ fix } \partial C)$ as a preimage, which is central.
The (1) and (2) prescribe a section of the whole group $\tilde A_f$. 

This forces $\pi_0 \Diff(C \text{ fix } \partial C)$ to be the semi-direct product of $\tilde A_f$
and $\prod_{i=1}^n \pi_0 \Diff(C_{J_i} \text{ fix } \partial C_{J_i})$.  The action
of $\tilde A_f$ on $\prod_{i=1}^n \pi_0 \Diff(C_{J_i} \text{ fix } \partial C_{J_i})$ 
(by conjugation) sends
$\Diff(C_{J_i} \text{ fix } \partial C_{J_i})$ to $\Diff(C_{J_k} \text{ fix } \partial C_{J_k})$
if the corresponding permutation $\tilde A_f \to A_f \to \Sigma_n^+$ sends $i$ to $k$.  If 
the permutation sends $i$ to $k$ preserving the sign, this means the $J_i$ and $J_k$ 
are isotopic as long knots and 
$\Diff(C_{J_i} \text{ fix } \partial C_{J_i}) \to \Diff(C_{J_k} \text{ fix } \partial C_{J_k})$
is conjugation by a diffeomorphism $C_{J_i} \to C_{J_k}$ that can be realized by an isotopy 
from the long knot $J_i$ to $J_k$. If on the other hand the sign is not preserved, then the map
$\Diff(C_{J_i} \text{ fix } \partial C_{J_i}) \to \Diff(C_{J_k} \text{ fix } \partial C_{J_k})$
is conjugation by a map $C_{J_i} \to C_{J_k}$ which can be realized by an isotopy in $\Real^3$
that reverses the long axis of the knot. Apply the classifying space functor $B$ to these
results and the theorem follows.
\end{proof}
\end{thm}

It's interesting to compare this result to a theorem that appears
in Hatcher's papers \cite{Hatcher4, Hatcher5} on
the homotopy type of the component of a hyperbolic knot $f$ in
$\Emb(S^1,S^3)/\Diff(S^1)$. Provided 
the Linearization Conjecture is true, one can put $f$ into
a `maximally-symmetric position,' which means that $f$ has a position
so that the restriction map
$\Isom_{S^3}(S^3,f) \to \Isom_{\Bbb H^3}(S^3 \setminus f)$ is an isomorphism.
$\Isom_{S^3}(S^3,f)$ is the group of isometries of $S^3$ that preserve
the knot, and $\Isom_{\Bbb H^3}(S^3 \setminus f)$ is the group of hyperbolic
isometries of $S^3 \setminus f$. Let $\Isom^+_{S^3}(S^3,f)$ be the subgroup 
of $\Isom_{S^3}(S^3,f)$ which preserves the orientation of $S^3$, 
thus $\Isom^+_{S^3}(S^3,f) \subset SO_4$. Hatcher's theorem is that the 
component of $f \in \Emb(S^1,S^3)/\Diff(S^1)$ has the 
homotopy type of 
$SO_4/\Isom^+_{\Bbb H^3}(S^3 \setminus f)$, where the representation
$ \Isom^+_{\Bbb H^3}(S^3\setminus f) = \Isom^+_{S^3}(S^3,f) \to SO_4$ is
given by the inclusion map. So not only does the homotopy type of the
components of $\Emb(S^1,S^3)/\Diff(S^1)$ depend on the isometry groups
$\Isom^+_{\Bbb H^3}(S^3 \setminus f)$, but also on the representation
$\Isom^+_{S^3}(S^3, f) \to SO_4$. 

Observe that the representation $B_L \to \Sigma_n^{+}$ is generally not faithful. 
We will give some examples, with $L=(L_0,L_1,\cdots,L_n)$ a hyperbolic KGL.
\begin{itemize}
\item For $n=0$ the figure-8 knot gives
$B_L = \Zed_2$ yet $B_L \to \Sigma_0^+$ is trivial. 
\item For $n=1$ the link called $6_3^2$ in Rolfsen's tables 
satisfies $B_L = \Zed_2$ with $B_L \to \Sigma_1^+$ trivial.
$$\includegraphics[width=3cm]{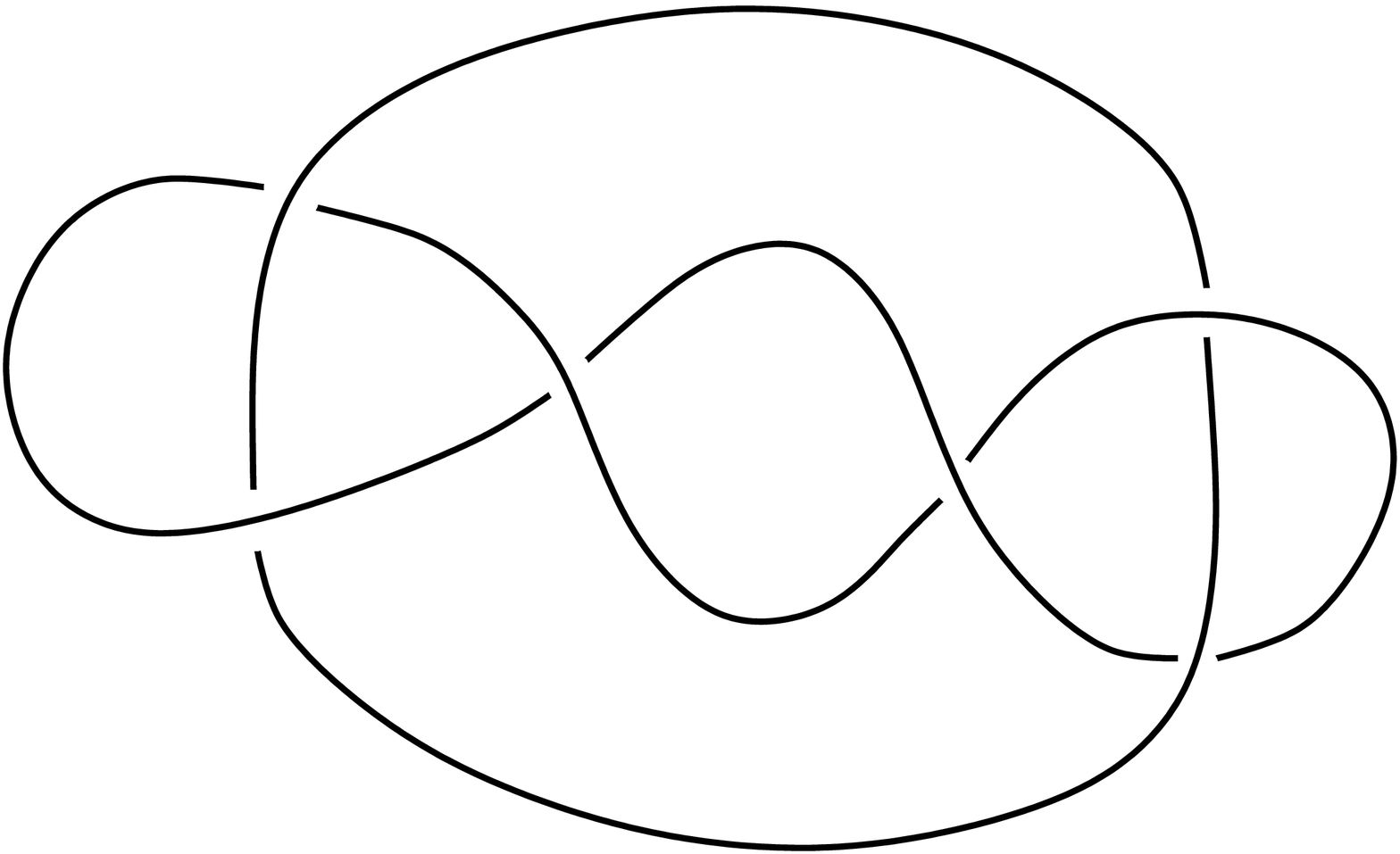}$$
\item For $n=2$ in the example below, $B_L=\Zed_4$ and
the kernel of $B_L \to \Sigma_2^+$ has order $2$
$$\includegraphics[width=3cm]{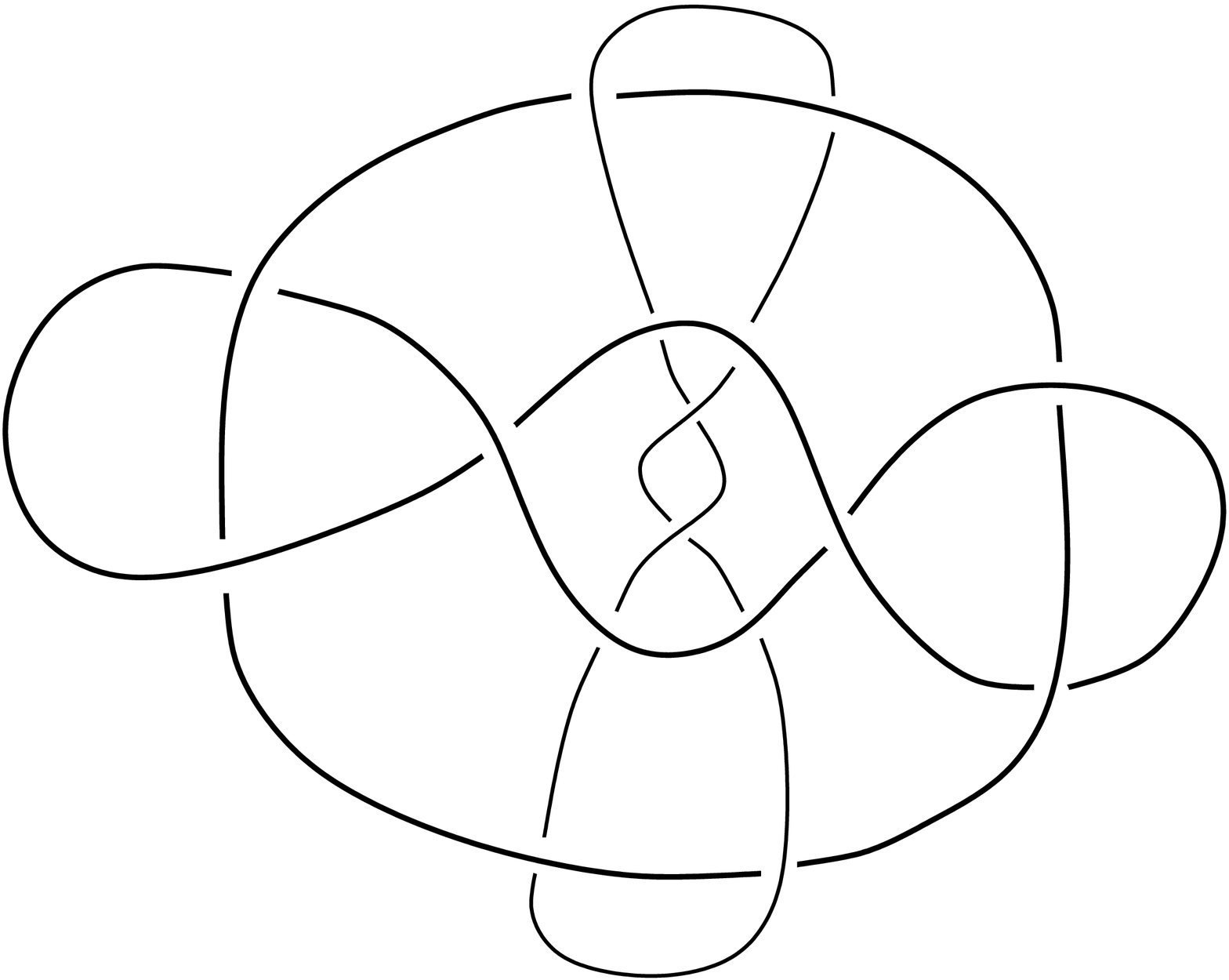}$$
\end{itemize}

Notice that in order to compute the homotopy type of $\K_f$,
if $A_f \to \Sigma_n^+$ is not conjugate to a representation
 $A_f \to \Sigma_n$, then we must know some of the spaces 
$\K_{J_i}$ as $O_2$-spaces. Minimally, we need to know
the homotopy class of the inversion map.  
In the next section we will inductively study the $O_2$
homotopy type of $\K$, and recursively compute the homotopy class of
the inversion map.

\section{Closed knots in $S^3$ and the $O_2$-action on $\K$}\label{3sphere}

This section describes the homotopy type of the space of
embeddings $\Emb(S^1,S^3)$.  This problem is reduced
to computing the homotopy type of $\K$ as an $O_2$-space, where
$O_2 \subset SO_3$ is the `long axis' preserving subgroup.

\begin{prop}\cite{BC}
$\Emb(S^1,S^3) \simeq SO_4 \times_{SO_2} \K$.
Moreover, there is a homotopy equivalence
$$SO_4 \times_{SO_2} \K \simeq S^3 \times \left( SO_3 \times_{SO_2} \K\right)$$
This makes $SO_3 \times_{SO_2} \K$ the $\K$-bundle over $S^2$ whose
monodromy is induced by the monodromy for the tangent bundle of $S^2$.
\end{prop}

Consequently, to describe the homotopy of the space of embeddings $\Emb(S^1,S^3)$
we need to know the homotopy type of $\K$ as an $SO_2$-space.  As we've seen, 
we can use Theorem \ref{mainthm} to compute the homotopy type of the components 
of $\K$ provided we know the homotopy type of simpler components of
$\K$, the $SO_2$-action on those components, and the homotopy class of the 
inversion involution $I : \K \to \K$ on those simpler components (provided the 
knot is invertible). So in this section we `kill two
birds with one stone,' and study $\K$ as an $O_2$-space.

\begin{prop}\label{eq1}{\bf ($O_2$-action on a cabling)}
If $f$ is a cable of $g$ then
\begin{itemize}
\item There is an $SO_2$-equivariant homotopy
equivalence

$$S^1 \times \K_g \to \K_f$$

Where the $SO_2$ action on $S^1 \times \K_g$ is given by $A\cdot (v,x)=(A\cdot v,x)$.

\item $I(\K_f)=\K_f$ if and only if
$I(\K_g)=\K_g$. If $I(\K_f)=\K_f$ then the above homotopy equivalence is 
$O_2$-equivariant where $A\cdot(v,x) = (A\cdot v,I^{\epsilon(A)}(x))$. Where
$\epsilon : O_2 \to \{\pm 1\}$ is the unique epimorphism.

\end{itemize}
\begin{proof}
Consider $f$ to be a splice knot $f=g \splice L$ where $L=S^{(p,q)}$ is a 2-component Seifert 
link \cite{jsj}. Theorem \ref{cable_thm} gives us a homotopy-equivalence
$SO_2 \times \K_g \to K_f$.  We will construct this map explicitly on the space-level. 
Consider the space of embeddings of $\Real \sqcup S^1$ in $\Real^3$ such that
the restriction to $\Real$ is a long knot.  Denote this space by $\Emb(\Real \sqcup S^1,\Real^3)$.
Notice that the Theorem \ref{fundprop} generalizes immediately to say that the component
of $\Emb(\Real \sqcup S^1,\Real^3)/\Diff(S^1)$ that contains the Seifert link $S^{(p,q)}$
has the homotopy-type of $S^1$.  Moreover, by the same argument as in Theorem \ref{cable_thm}, 
this component has the $SO_2$-equivariant homotopy-type of $S^1$ with the standard $SO_2$ action.
Since $f = g \splice S^{(p,q)}$, we can use parametrized splicing \cite{jsj} to construct an 
$SO_2$-equivariant map $S^1 \times \K_g \to \K_f$ which is also by design a homotopy-equivalence. 
\end{proof}
\end{prop}

Proposition \ref{eq1} also applies to a torus knot component.  In this
case, it gives the result that a torus knot component has the homotopy type of
$S^1$ and the $O_2$-action on $S^1$ is precisely the standard action of
$O_2$ on $S^1$.  The next proposition assumes familiarity with the action 
of the operad of $2$-cubes on $\K$ \cite{cubes}.

\begin{prop}\label{eq2}{\bf ($O_2$-action on a connected-sum)}
\begin{itemize}

\item By \cite{cubes}, if $f = f_1 \# \cdots \# f_n$ where 
$\{f_i : i\in \{1,2,\cdots,n\}\}$ are the prime summands of $f$, then there
is a canonical homotopy equivalence
$$\Cu_2(n) \times_{\Sigma_f} \prod_{i=1}^n
\K_{f_i} \to \K_f$$
This homotopy equivalence is
$SO_2$-equivariant, where the $SO_2$ action on $\Cu_2(n) \times_{\Sigma_f}
\prod_{i=1}^n \K_{f_i}$ is given by
$$A\cdot((C_1,\cdots,C_n),(g_1,\cdots,g_n))=((C_1,\cdots,C_n),(A\cdot g_1,\cdots,A\cdot g_n))$$

\item  The $O_2$ action on $\K$ restricts to $\K_f$
if and only if there is some involution $\iota \in \Sigma_n$ so that
$\K_{If_i} = \K_{f_{\iota(i)}}$ for all $i \in
\{1,2,\cdots,n\}$. 
In this case the map 
$\Cu_2(n) \times_{\Sigma_f} \prod_{i=1}^n \K_{f_i} \to \K_{f}$ is
$O_2$-equivariant where the $O_2$ action on 
$\Cu_2(n) \times_{\Sigma_f} \prod_{i=1}^n \K_{f_i}$ is given by
$$A\cdot ((C_1,\cdots,C_n),(g_1,\cdots,g_n))=
 (\epsilon(A) \cdot (C_{\iota(1)},\cdots,C_{\iota(n)}),(A\cdot g_{\iota(1)},\cdots,A \cdot g_{\iota(n)}))$$
where $\epsilon(A) \in \{\pm 1\}$ acts on $\Cu_2(n)$ by mirror reflection
about the $y$-axis and $\epsilon : O_2 \to \{\pm 1\}$ is the 
unique epimorphism. 
\end{itemize}
\begin{proof}
The homotopy equivalence is the little $2$-cubes action $\cubeact$ \cite{cubes}.
To check equivariance, we use the 
`trivially framed' space of knots $\tK$ \cite{cubes} which admits the 
$2$-cubes action $\cubeact$ and is homotopy equivalent to $\K$.
First, if $A \in SO_2$, $(L_1,L_2,\cdots,L_n)\in \Cu_2(n)$, $(f_1,f_2,\cdots,f_n)\in \tK^n$ then
$$ A\cdot\cubeact((C_1,C_2,\cdots,C_n),(g_1,g_2,\cdots,g_n)) 
=  (Id_{\Real} \times A) \circ (C_{\sigma(1)}^\pi\cdot g_{\sigma(1)} \circ \cdots \circ C_{\sigma(n)}^\pi\cdot g_{\sigma(n)}) \circ (Id_{\Real} \times A^{-1})$$
$$=  (Id_{\Real} \times A)  \circ C_{\sigma(1)}^\pi\cdot g_{\sigma(1)} \circ 
(Id_{\Real} \times A^{-1}) \circ \cdots \circ (Id_{\Real} \times A)\circ 
C_{\sigma(n)}^\pi\cdot g_{\sigma(n)} \circ (Id_{\Real} \times A^{-1})$$
$$ = C^{\pi}_{\sigma(1)}\cdot ( (Id_{\Real}\times A)\circ g_{\sigma(1)} \circ (Id_{\Real}\times A^{-1}))\circ \cdots \circ  C^{\pi}_{\sigma(n)}.
( (Id_{\Real}\times A)\circ g_{\sigma(n)} \circ (Id_{\Real}\times A^{-1})) $$
$$ = C^{\pi}_{\sigma(1)}\cdot ( A\cdot g_{\sigma(1)}) \circ \cdots \circ  C^{\pi}_{\sigma(n)}.
( A\cdot g_{\sigma(n)})=\cubeact((C_1,C_2,\cdots,C_n),A\cdot (g_1,g_2,\cdots,g_n))$$
In the above sequence of equations, the 3rd equality is perhaps the only step where
a little explanation is required.
$$(Id_\Real \times A)\circ C^\pi_{\sigma(i)}\cdot g_{\sigma(i)} \circ (Id_\Real \times A^{-1})
= C^\pi_{\sigma(i)}\cdot ((Id_\Real \times A)\circ g_{\sigma(i)} \circ (Id_\Real \times A^{-1}))$$
Since both $C^\pi_{\sigma(i)}$ and $Id_{\Real}\times A$ act on the knot
space as a product action, $C^\pi_{\sigma(i)}$ acting only in the direction of long axis
(identity in the orthogonal directions), while $Id_{\Real}\times A$ is by design the identity
on the long axis, they commute.

To justify the equivariance with respect to inversion is essentially the same 
computation. If $C_i$ is a 2-cube, denote $C_i$ conjugated by the map
$(x,y) \longmapsto (-x,y)$ by $C'_i$. 

$$I\cdot \cubeact((C_1,C_2,\cdots,C_n),(g_1,g_2,\cdots,g_n)) = \cubeact( (C'_1,\cdots,C'_n), (I\cdot g_1,\cdots, I\cdot g_n))$$
\hskip 20mm $=\cubeact( (C'_{\iota(1)}, \cdots,C'_{\iota(n)}), (I\cdot g_{\iota(1)},\cdots,I\cdot g_{\iota(n)}))$.

The first step is argued similarly to the previous case of the $SO_2$-action, the second 
step is the symmetry axiom of a cubes action.
\end{proof}
\end{prop}

Note in Proposition \ref{eq2}, the choice of involution $\sigma$ is irrelevant, since
any two choices differ by an element of $\Sigma_f$.

\begin{prop}\label{eq3}{\bf ($O_2$-action on a hyperbolically spliced knot)}
If $f=J \splice L$ where $L$ is a hyperbolic KGL:
\begin{itemize}
\item The homotopy equivalence given in Theorem \ref{mainthm}
$$SO_2 \times \left( (\prod_{i=1}^n \K_{J_i})   \times_{A_f} L_0 \right) \to \K_f$$
can be made to be $SO_2$-equivariant, where the $SO_2$ action on
$SO_2 \times \left( (\prod_{i=1}^n \K_{J_i}) \times_{A_f} L_0 \right)$ is a product action, 
$SO_2 \times SO_2 \to SO_2$ is left multiplication and the $SO_2$-action on 
$(\prod_{i=1}^n \K_{J_i}) \times_{A_f} L_0$ fixes all points.
\item $I(\K_f)=\K_f$ if and only if there is an isotopy from $L$ to $I(L)$ as unoriented
links satisfying: 
\begin{itemize}
\item The isotopy is orientation-preserving between  $L_0$ and $I(L_0)$.
\item Given the above isotopy, define a signed permutation $\iota \in \Sigma_n^{+}$ 
so that $\iota(i)= \epsilon j$ where $\epsilon \in \{\pm 1\}$
provided the isotopy sends $L_i$ to $\epsilon(i)\cdot L_j$ where here $1\cdot L_j=L_j$ and
$-1\cdot L_j$ is $L_j$ with orientation reversed. 
\item $\iota\cdot (\K_{J_1}\times\K_{J_2}\times \cdots \times \K_{J_n})= \K_{J_1}\times\K_{J_2}\times \cdots \times \K_{J_n}$.
\end{itemize}

In this case, $I$ is conjugate (via our homotopy-equivalence) to
the self-map of $SO_2 \times \left( (\prod_{i=1}^n \K_{J_i}) \times_{A_f} L_0 \right)$
given by: 
$$(A, (J_1,J_2,\cdots,J_n), x) \longmapsto (c(A), \iota\cdot (J_1,J_2,\cdots,J_n), i(x))$$
where $c : SO_2 \to SO_2$ and $i : L_0 \to L_0$ are any orientation-reversing involutions. 
 \end{itemize}
\begin{proof}
A knot is invertible if and only if $I(\K_f)=\K_f$, or equivalently,
if and only if $C$ admits an orientation-preserving diffeomorphism $g : C \to C$
which reverses the longitudinal orientation of $f$. Such a
diffeomorphism can be isotoped to preserve the JSJ-decomposition of $C$. Thus $g$ restricts 
to a diffeomorphism of $V$ and a diffeomorphism of $\sqcup_{i=1}^n C_{J_i}$, 
defining a signed permutation $\iota \in \Sigma_n^+$ by how it permutes components and longitudinal
orientations of $\sqcup_{i=1}^n C_{J_i}$.   The automorphism of $V$ acts as an involution 
(by conjugation) on $\tilde A_f$. This action is negation, since it restricts to negation on the meridional
Dehn Twist.  

The homotopy equivalence
$$SO_2 \times \left( (\prod_{i=1}^n \K_{J_i}) \times_{A_f} L_0 \right) \to \K_f$$
can be made $SO_2$-equivariant via a rather small observation:  Restrict
the above homotopy equivalence to a map a point in $SO_2$, this gives a map
$F: (\prod_{i=1}^n \K_{J_i}) \times_{A_f} L_0 \to \K_f$.  Since $SO_2$ acts
on $\K_f$ we get a new map 
$SO_2 \times \left( (\prod_{i=1}^n \K_{J_i}) \times_{A_f} L_0 \right) \to \K_f$ by
sending $(A,x)$ to $A\cdot F(x)$.  This map is also a homotopy-equivalence, 
since $\tilde A_f$ splits as a direct sum of two infinite cyclic groups, 
with one of the generators a meridional Dehn twist. 
\end{proof}
\end{prop}

It seems likely that in the case of an invertible knot, 
the homotopy equivalence in Proposition \ref{eq3} could be made $O_2$-equivariant for a suitable
action of $O_2$ on  $SO_2 \times (\prod_{i=1}^n \K_{J_i}) \times_{A_f} L_0$, 
although a few computations with rather simple knots such as the Whitehead double
of a non-invertible hyperbolic knot demonstrates that with the specific 
action given in Proposition \ref{eq3} does not suffice \cite{Hartley}. 
Finding an $O_2$-equivariant homotopy equivalence is likely a subtle
problem.

\section{Examples}\label{Examples}

How difficult is it to compute the homotopy type of a component of $\K$?  
The answer depends on what your starting point is.  
The simplest starting point would be a companionship-tree description 
of a knot $f$, together with a list of the relevant hyperbolic
isometry groups $B_L$ and their representations $B_L \to \Sigma^+_n$.
In that case, one can simply read off the homotopy type of $\K_f$ from
the theorems in this paper. For the sake of completeness, one also 
requires knowledge of the groups $\Isom_{\Bbb H^3}(S^3 \setminus L)$ 
to determine invertibility and the homotopy class of the inversion map
of the various companion knots to $f$.

If on the other hand one takes the more traditional starting point of
a Reidemeister knot-diagram for $f \in \K$, one needs to first compute 
the companionship tree $G_f$.  The algorithm of Jaco, Letscher, 
Rubinstein \cite{JLR} allows one to compute $G_f$.
The JLR algorithm recognises the Seifert-fibred pieces together
with the singular fibre data. The remaining pieces are hyperbolisable. 
In practise, one frequently can compute the hyperbolic structure
and the representation $B_L \to \Sigma_n^+$ with SnapPea.  As of yet,
there is no proof that a SnapPea-type procedure \cite{Snap}
will find the hyperbolic structure on a hyperbolisable manifold.
There is an algorithm to find the hyperbolic structure.
Haken and Hemion have a solution to the word problem in
Haken manifolds \cite{HakHem} and
Manning has an algorithm to find the hyperbolic structure
on a hyperbolisable manifold, provided it comes equipped
with a solution to its word problem \cite{Man}. In practise this
approach seems rather difficult to implement.

By and large, the computation of the representations $B_L \to \Sigma_n^+$
is most easily handled by SnapPea. It searches for and frequently
finds the canonical polyhedral 
decomposition of a cusped hyperbolic manifold $M$ \cite{Epstein}.
The group of hyperbolic isometries of $M$ is precisely the group
of combinatorial automorphisms of the canonical polyhedral 
decomposition. 

If one wishes to avoid the use of computer programs, Kouzi
and Sakuma show how one can compute the symmetry 
groups of various families of hyperbolic links in $S^3$ \cite{Sakuma2}  
by using the $\Zed_2$-equivariant JSJ-decomposition of the $2$-sheeted 
branched cover of $S^3$, branched over the link. 
This is frequently quite practical.

If one takes yet a different starting point, it's potentially much 
easier to `compute' the homotopy type of $\K$.  At present the
class of representations $B_L \to \Sigma_n^+$ for hyperbolic KGLs
is not understood.  But if one had a list of all such representations,
the homotopy type of $\K$ could be described without
ever mentioning a knot or a link -- it would simply
be the recursive closure of the total spaces in 
Propositions \ref{eq1}, \ref{eq2} and \ref{eq3}, generated
by a point space and our prescribed class of representations
$B_L \to \Sigma_n^+$.

We take the first point of view in this section, and explicitly 
compute the homotopy type of $\K_f$. 

{\bf Some examples:}

\begin{itemize}
\item $\K_f \simeq \{*\}$ if $f$ is the unknot.
\item $\K_f \simeq S^1$ if $f$ is a non-trivial torus knot. Here
$S^1$ is represented by the rotations of a fixed torus knot about the long 
axis.

{
\psfrag{torusknot}[tl][tl][0.8][0]{\hskip -6mm A torus knot $f$, with $\K_f \simeq S^1$}
$$\includegraphics[width=4.5cm]{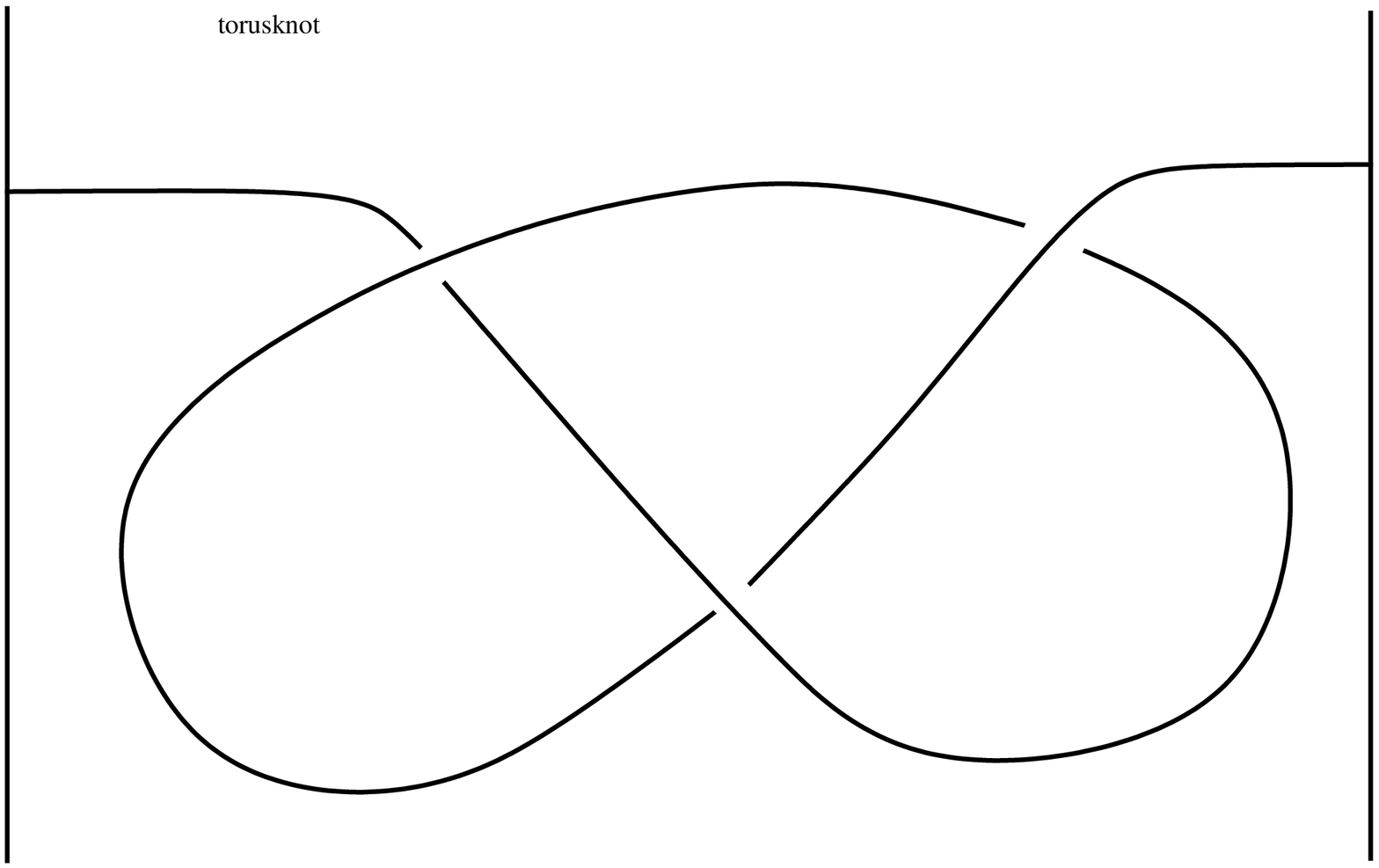}$$
}

\item $\K_f \simeq (S^1)^n$ if $f$ is an $n$-fold iterated cable of the unknot
ie: an $(n-1)$-fold cable of a torus knot. In the figure below, $\K_f \simeq (S^1)^2$
is represented by the rotations of the knot about the long axis, and the rotation of the 
cabling parameter respectively ($\K_g$ in Proposition \ref{eq1}).
{
\psfrag{cableknot}[tl][tl][0.9][0]{\hskip -2mm $\K_f \simeq (S^1)^2$}
$$\includegraphics[width=4cm]{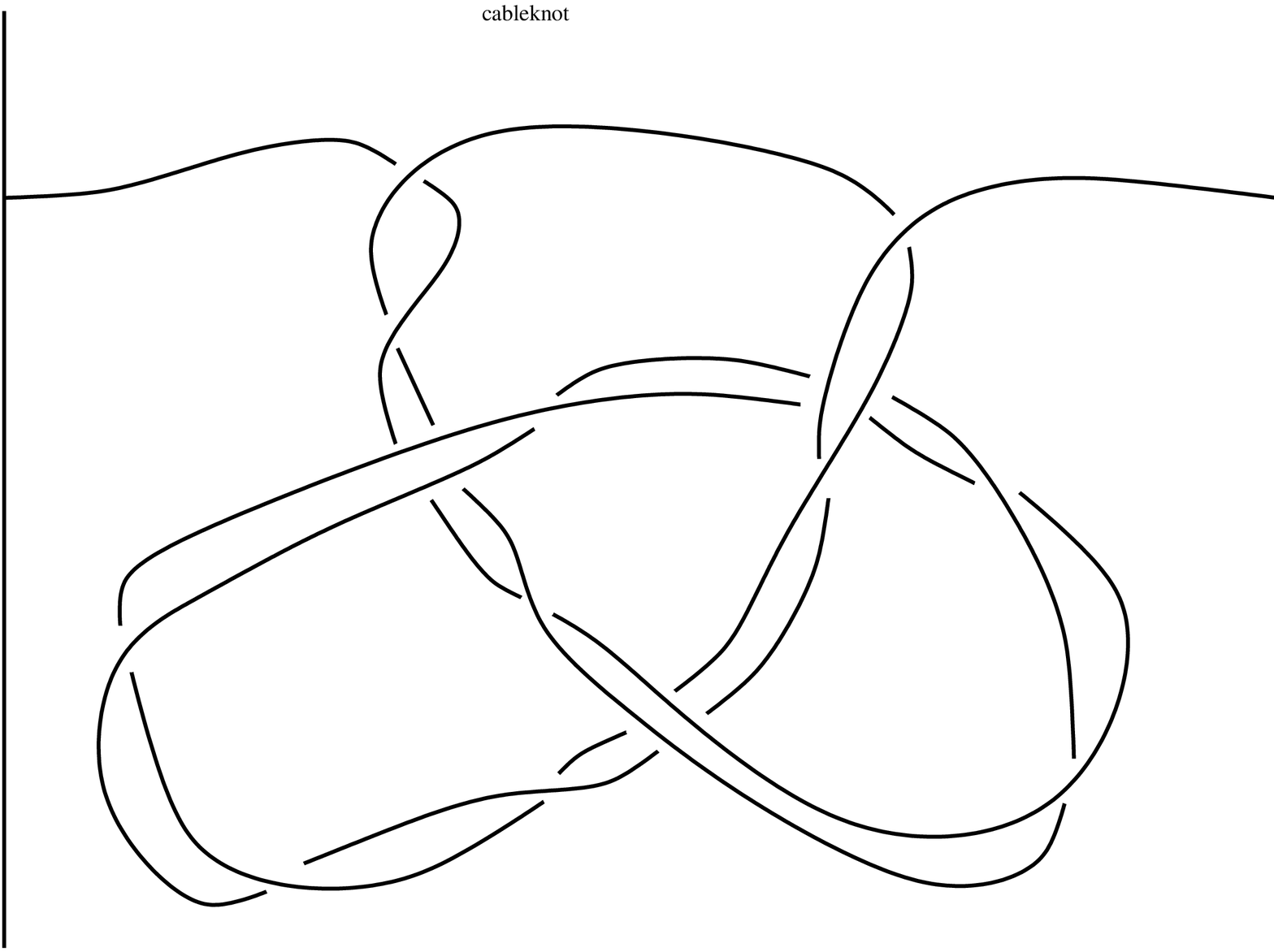}$$
}
\end{itemize}

A few of the remaining examples use a signed cycle notation for elements
of $\Sigma_n^+$.
The signed cycle $(a_1 a_2 \cdots a_k \pm) \in \Sigma_n^+$ is the signed
permutation where $a_1 \longmapsto a_2$, $a_2 \longmapsto a_3$, preserving
signs, but $a_k \longmapsto \pm a_1$, meaning the sign is preserved if
$\pm = +$ and reversed if $\pm = -$. 

For the Borromean rings $L$, $B_L \simeq \Zed_4$
where the representation $B_L \to \Sigma_2^+$ sends the generator
to the signed 2-cycle $(12-)$, which has order 4.   To those that are
curious, the full group of isometries of the complement of the Borromean
rings has order $48$.  This includes orientation-reversing isometries.  The
subgroup that preserves the orientation of $S^3 \setminus L$ has order
$24$.  There is an element of order $3$ which cyclically permutes the 
components of $L$, which can be seen in the `standard' picture of $L$.
$$\includegraphics[width=3cm]{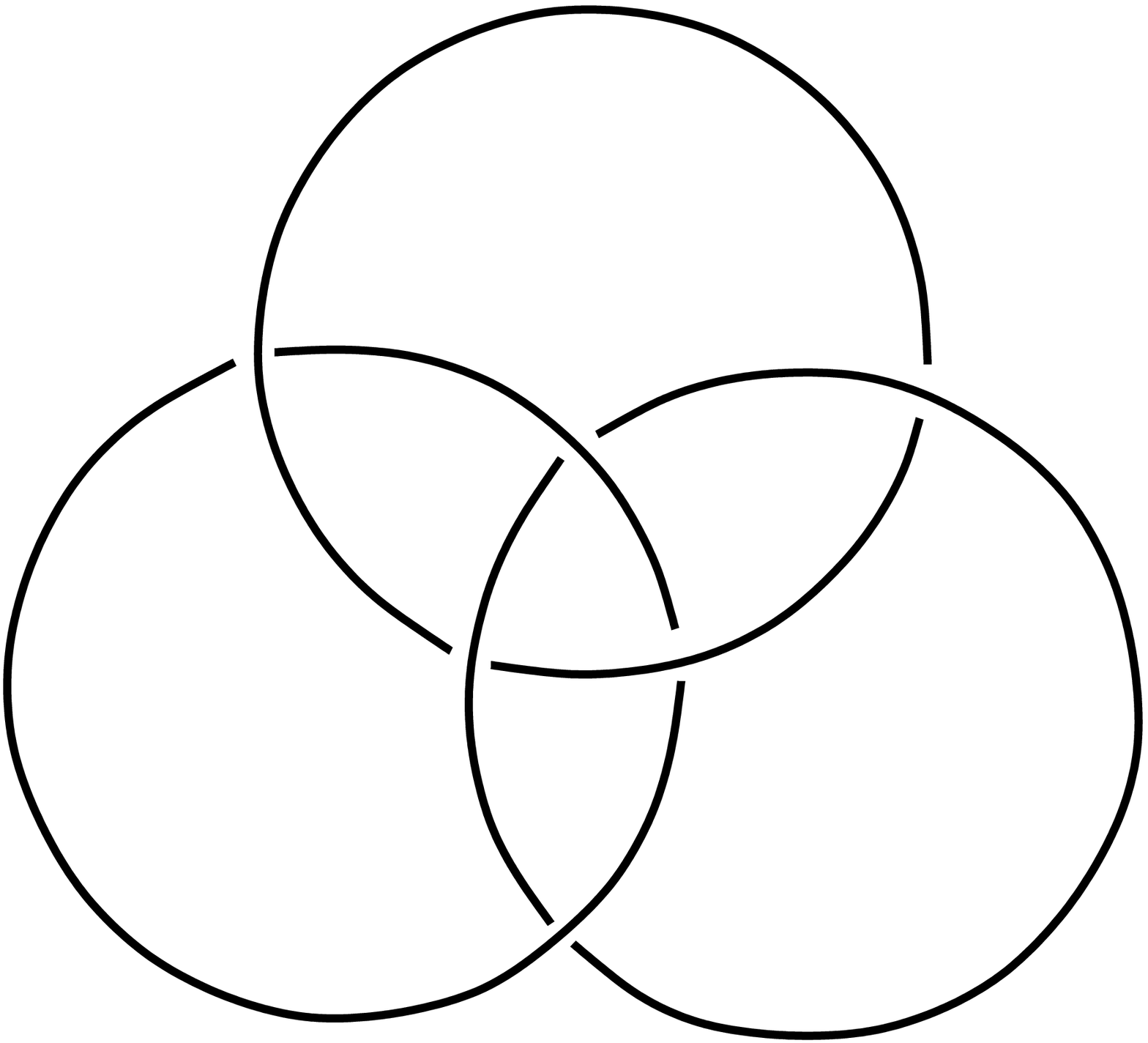}$$
There is an element of order $2$ which preserves a component of $L$ but
reverses its orientation.  This leaves $B_L$ as a subgroup of order $4$. 
This symmetry is perhaps the most difficult to visualize of the symmetries 
of $S^3 \setminus L$, though it and the full symmetry group of order $48$
can be realized as isometries of $S^3$ when $L$ is put in the `atomic symbol'
position.
$$\includegraphics[width=3cm]{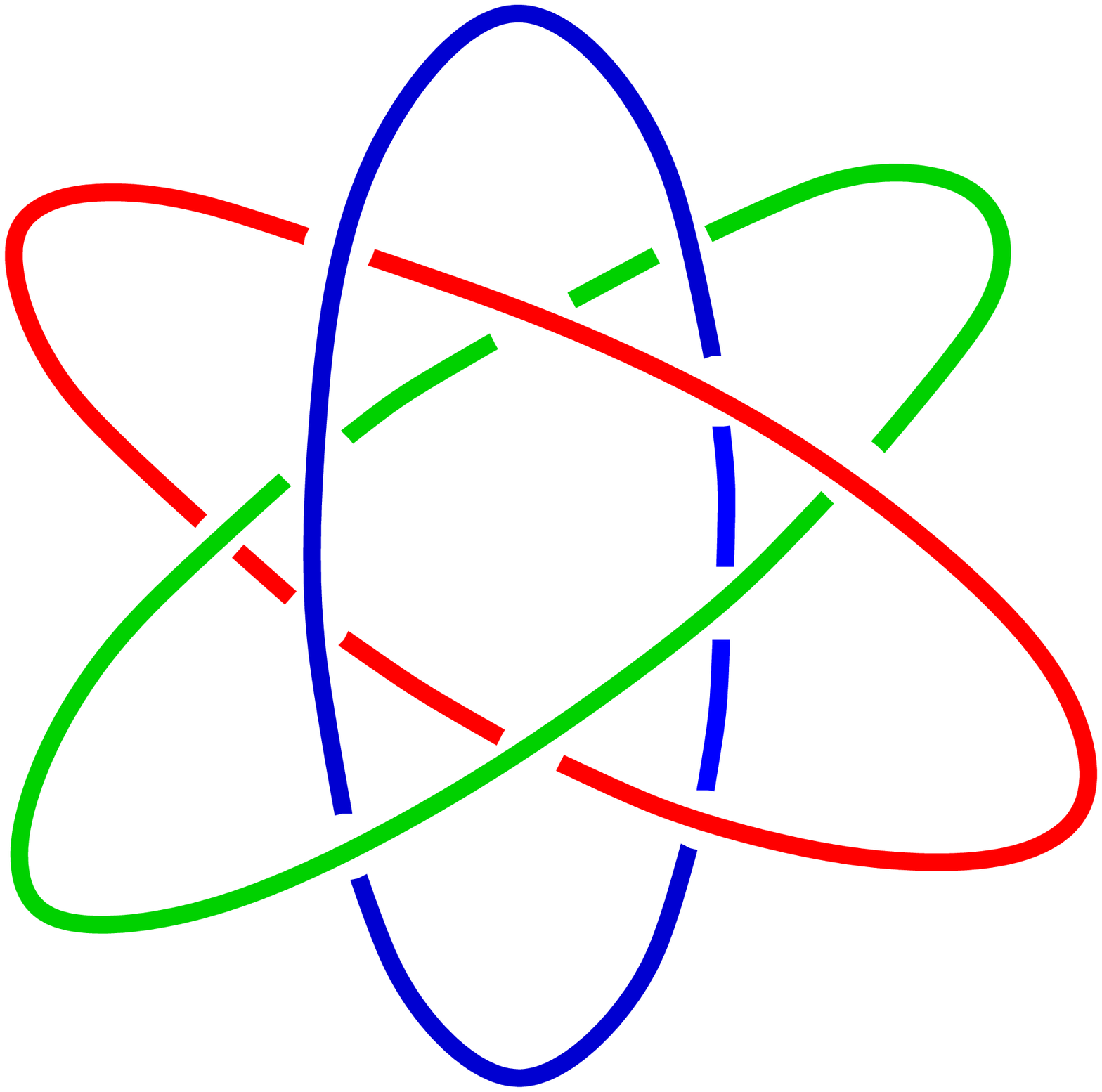}$$

\begin{itemize}
\item Consider a `Borromean sum' of two figure-8 knots,
$f = (J_1,J_2) \splice L$  
where $L$ is the Borromean rings, and $(J_1,J_2)$ are two copies of
the figure-8 knot then $A_f = B_L$ and
$\K_f \simeq S^1 \times (S^1 \times_{\Zed_4} ((S^1)^2 \times (S^1)^2))$.
 The action of the generator of $\Zed_4$ on the fibre
is induced by the matrix 
$\left( \begin{matrix} 0 & -1 \cr 1 & 0 \end{matrix} \right)$, where
$-1$ acts on $S^1 \times S^1$ by rotation by $\pi$ about some
fixed point. Thus, $\K_f$ is a $6$-dimensional
Euclidean manifold.
{
\psfrag{cableknot}[tl][tl][0.8][0]{A cable of a torus knot}
$$\includegraphics[width=7cm]{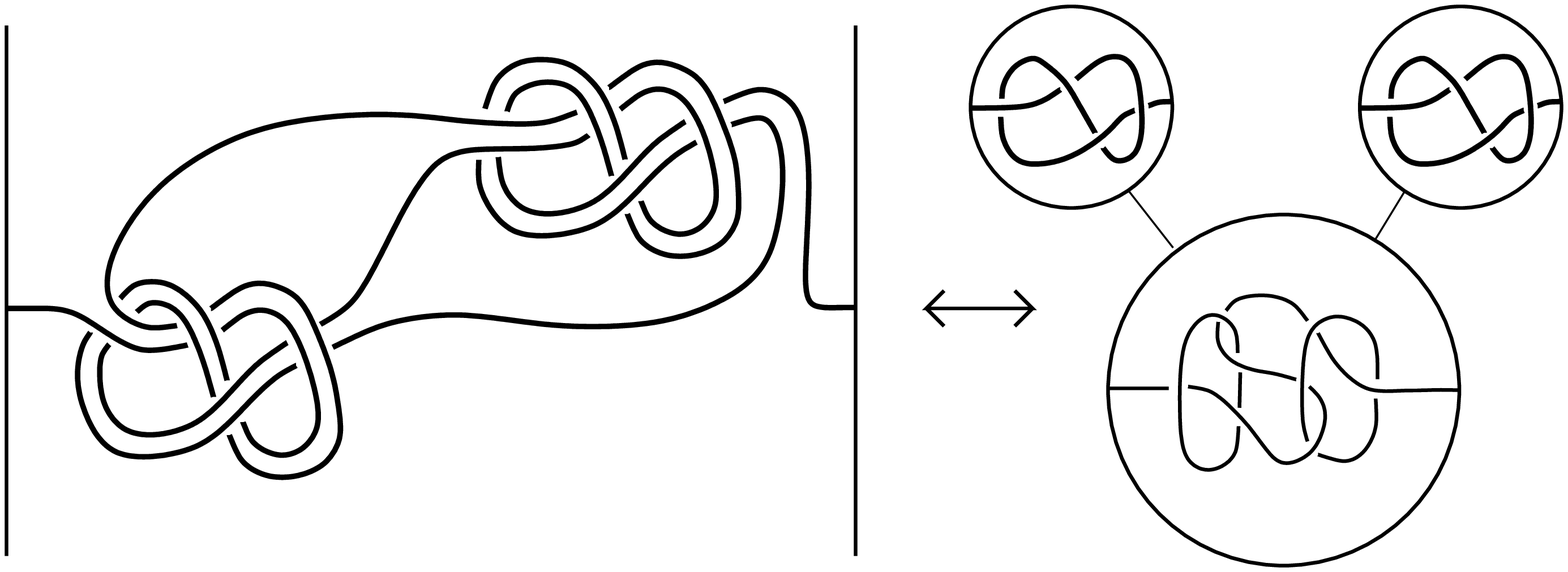}$$
}
\item In the above case, if we replace one of the figure-8 knots
by a trefoil we get $A_f \simeq \Zed_2$ so 
$\K_f \simeq S^1 \times (S^1 \times_{\Zed_2} (S^1 \times (S^1)^2))$.
$\Zed_2$ acts on $S^1$ by conjugation, on $(S^1)^2$ by rotation by $\pi$
about some fixed point.
{
\psfrag{f}[tl][tl][0.8][0]{}
\psfrag{J1}[tl][tl][0.8][0]{}
\psfrag{J2}[tl][tl][0.8][0]{}
\psfrag{L}[tl][tl][0.8][0]{}
$$\includegraphics[width=7cm]{borromean.eps}$$
}
\end{itemize}
Consider a Whitehead link $L$, then $B_L = \Zed_2$.
If $f$ is obtained from $J$ by Whitehead doubling, ie: 
$f=J \splice L$, then:

\begin{itemize}
\item $\K_f \simeq S^1 \times (S^1 \times_{\Zed_2} \K_{J})$ if $\K_{J}=\K_{IJ}$
{
\psfrag{n}[tl][tl][0.8][0]{}
\psfrag{GHL}[tl][tl][0.8][0]{}
$$\includegraphics[width=6cm]{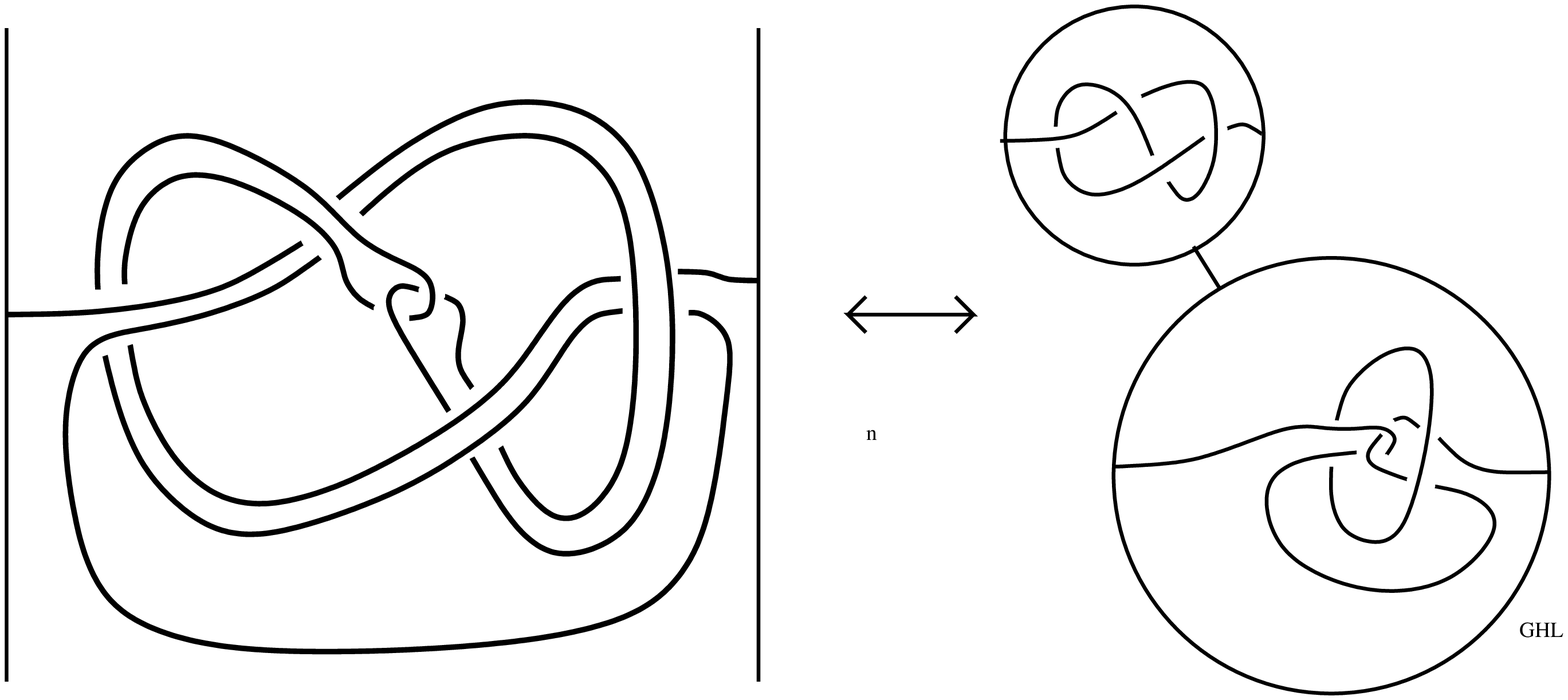}$$
}
\item If $J$ is not invertible, then $A_f$ is trivial and
 $\K_f \simeq S^1 \times S^1 \times \K_J$.
\end{itemize}

For the knot below, it's a good exercise to compute the 
homotopy type of $\K_f$.
{
$$\includegraphics[width=5cm]{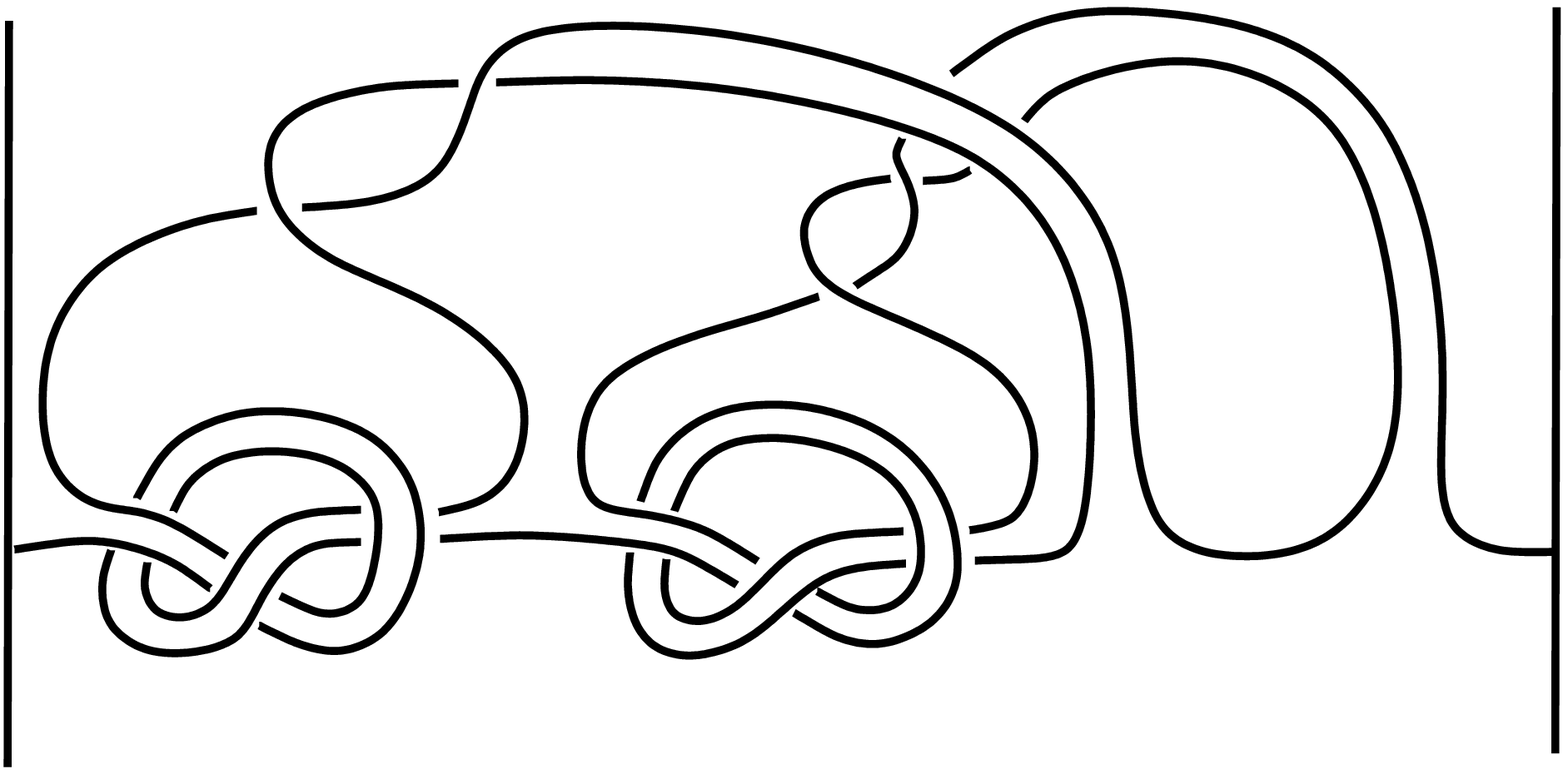}$$
}

\begin{itemize}
\item A more complicated example with interesting homological behavior is given by letting 
$f_1 \in \K$ be a trefoil knot. Then $\K_{f_1}$ has the
homotopy type of a circle. Let $f_2$ be the 
connect-sum of $4$ copies of $f_1$, $f_2= \#_4 f_1 \in \K$. 
Then $\K_{f_2}$ has the homotopy type of $\Cu_2(4) \times_{\Sigma_4} (S^1)^4$.
Let $f_3$ be a cable of $f_2$, then $\K_{f_3}$ has the homotopy type of
$S^1 \times \left( \Cu_2(4) \times_{\Sigma_4} (S^1)^4 \right)$. Let $f_4$
be the connect-sum of $2^n$ copies of $f_3$, then $\K_{f_4}$ has the homotopy type
of $\Cu_2(2^n) \times_{\Sigma_{2^n}} \left( S^1 \times \left( \Cu_2(4) \times_{\Sigma_4} (S^1)^4 \right)\right)^{2^n}$. In \cite{BC} we prove that the homology of $\K_{f_4}$
has torsion of order $2^{n+1}$. Moreover, we show for any prime power $p^n$, one
can find a component of $\K$ whose homology contains torsion of order $p^n$.
\end{itemize}

\section{On the Gramain subgroup}\label{gramelt}

In Gramain's paper \cite{Gramain2} he points out that $\pi_1 \K_f$ contains a group
isomorphic to the integers if $f$ is a non-trivial knot.  This section investigates
that subgroup. 

Gramain's subgroup of $\pi_1 \K_f$ is given by the $SO_2$-action on $\K_f$.  Consider
the map $SO_2 \times \K_f \to \K_f$ and restrict it to $SO_2 \times \{f\}$.  The induced
map $\pi_1 SO_2 \to \pi_1 \K_f$ is the Gramain homomorphism, and the image is Gramain's 
subgroup.  Gramain proved that this map is an embedding provided $f$ is not the unknot. 
His argument was an application of the Burde-Zieschang characterization
of the centres of knot groups. We denote the Gramain subgroup of $\pi_1 \K_f$ by $G_f$. 
Sometimes we call generators of this group `the Gramain element.'

\begin{prop}\label{gramain_cocycle}
Provided $f$ is a non-trivial knot, the induced map $G_f \to H_1 \K_f$ is an embedding. 
Moreover, there is a map $\K_f \to SO_2$ such that the composite
$SO_2 \to \K_f \to SO_2$ has non-zero degree.
\begin{proof}
Since we have an inductive description of the homotopy type of the space of long knots, 
we prove this proposition by induction on the complexity of the JSJ-tree of the knot
complement.  Our base case is a tree with only one vertex -- this corresponds to either
a torus knot or a hyperbolic knot, in which case $\K_f \simeq S^1$ or 
$\K_f \simeq SO_2 \times \left(L_0/A_f\right)$. In the case of
the torus knot our map is simply the homotopy equivalence $\K_f \to S^1$. 
In the hyperbolic satellite case, our map is projection onto the meridional group
$\K_f \to SO_2$. 

We have three inductive steps, 
\begin{itemize}
\item $f$ is a cable knot.
\item $f$ is a connect-sum.
\item $f$ is a hyperbolic satellite.
\end{itemize}

In the case that $f$ is a cable knot, $\K_f \simeq S^1 \times \K_g$. The map
out of $\K_f$ is projection onto the $S^1$ factor.

In the case that $f$ is a connected-sum, 
$\K_f \simeq \Cu_2(n) \times_{\Sigma_f} \prod_{i=1}^n \K_{f_i}$. 
The map out of $\Cu_2(n) \times_{\Sigma_f} \prod_{i=1}^n \K_{f_i}$ 
is given by the product of all the maps $\K_{f_i} \to S^1$ where we
consider $S^1$ to be an abelian group.

In the case that $f$ is a hyperbolic satellite knot, we have a fibration 
a copy of $SO_2$ splits off $\K_f$, so we project onto $SO_2 \equiv S^1$.
\end{proof} 
\end{prop}

\begin{defn}
The image of the generator of $G_f$ in $H_1 \K_f$ is called the Gramain cycle, and 
will be denoted $g_*(f)$. The element of $H^1 \K$ constructed in Proposition
\ref{gramain_cocycle} is called the Gramain co-cycle and is denoted $g^*$. 
\end{defn}

\begin{prop}\label{gramobs} $\langle g^*,g_*(f)\rangle =
\left\{
\begin{array}{ll}
0 & \text{if } f \text{ is the unknot } \\
n & \text{if } f \text{ is a connect sum of } n \text{ prime knots } \\
\end{array}
\right.
$
\end{prop}

\section{Questions}\label{endgame}

In a paper or Tourtchine \cite{Turchin}, the `first non-trivial $1$-co-cycle'
detected by the Vassiliev spectral sequence for $\K$,  $v^1_3 \in H^1 \K$
satisfies:
\begin{itemize}
\item $\langle v^1_3,g_*(f)\rangle = v_2(f)$ mod $2$. 
\item $\langle v^1_3,\{a,b\}\rangle = 0$ mod $2$, where $\{\cdot,\cdot\}$ is the Poisson bracket on $H_* \K$ coming from the cubes action \cite{BC}.
\end{itemize}

$g^*$ is also trivial on Poisson brackets, so perhaps
$v_2\cdot g^*$ and $v_3^1$ agree on prime knot components?
By Proposition \ref{gramobs}, they generally do not agree
on connect-sums. $v_3^1$ appears rather mysterious from the point
of view of the computations contained in this paper.

An interesting direction to pursue would be to compute the
maps $\K \to AM_n$ for $n \geq 3$ coming from the Goodwillie Calculus
\cite{c2, Dev}. As pointed out by Goodwillie and Sinha \cite{Dev2}, the 
map $\K \to AM_2$ is null-homotopic. Since the components of $\K$ are $K(\pi,1)$'s, 
it would be interesting to compute the maps $\pi_1 \K_f \to \pi_1 AM_n$.
$\pi_1 AM_3 \simeq \Zed^2$, so it's possible that $\pi_1 AM_3$ detects the
Gramain element and perhaps one more element, such as projection onto
$T_l/A_f$ in the hyperbolic satellite case.  The map $\K \to \pi_0 AM_3 \simeq \Zed$
is the universal type-2 invariant of knots \cite{c2}, so perhaps
$\pi_1 AM_3$ fits better with Tourtchine's $v_3^1$ cocycle. 
 
One can give $\K$ the structure of a topological monoid in at least
two different ways. One method would be to think of $\K$ as a subspace
of $\EK{1,D^2}$ as in \cite{cubes}, then the monoid structure is 
induced by composition of functions.  Another way to endow $\K$ with
the structure of a topological monoid is to take the `Moore loop-space
model' for $\K$. Either way, the group completion 
$\Omega B \K$ of $\K$ is defined \cite{May, Paolo}
and has the homotopy type \cite{cubes, BC} of:
$$\Omega^2 \Sigma^2 (\Prime \sqcup \{*\}) \simeq 
  \Omega^2 \bigvee_{[f] \in \pi_0 \Prime} \left( S^2 \vee \Sigma^2 \K_f \right) $$ 
Since one can construct infinitely
many components of $\Prime$ that have the homotopy type of a 
product of circles (eg: iterated cable knots), 
we can conclude that $\Omega B\K$ has
$\Omega^2 \left( \bigvee_{i=2}^\infty \bigvee_{\infty} S^i \right)$ as a retract. 
By the results in \cite{BC}, $\Omega B\K$ does not have the same homotopy type
as $\Omega^2 \left( \bigvee_{i=2}^\infty \bigvee_{\infty} S^i \right)$ since
$H_*(\Omega B\K;\Zed)$ has torsion of all orders \cite{BC} while 
$H_* \left( \Omega^2 \left( \bigvee_{i=2}^\infty \bigvee_{\infty} S^i \right);\Zed\right)$
does not \cite{CLM}.

It may be reasonable to expect that the map
$\K \to AM_n$ factors as in the diagram below.
$$\xymatrix{\K \ar[r] \ar[dr] & \Omega B\K \ar[d] \\
                              & AM_n }$$
If so, this could be a good platform for studying convergence issues
for the Vassiliev spectral sequence, as the group completion map
$\K \to \Omega B\K$ induces an embedding on homology \cite{CLM}, while 
the fibres of $AM_n \to AM_{n-1}$ are iterated loop spaces on
total homotopy fibres of cubical diagrams of Fadell-Neuwirth fibrations \cite{Good}. 

Perhaps the most basic open question related to the homotopy type
of $\K$ is the issue of realization of representations $B_L \to \Sigma_n^+$
for hyperbolic KGLs $L$. We have seen a few examples where $B_L= \Zed_2$ or $\Zed_4$ and
a few representations of these groups, but if one spends a little more time
considering hyperbolic KGLs, it's clear this is far from all such representations.
Alexander Stoimenow has an interesting example
$B_L = \Zed_n \to \Sigma_n$ with generator mapping to an $n$-cycle.
All the links in the family are hyperbolic
since they are alternating, non-split and prime \cite{Menasco}.
{
\psfrag{n}[tl][tl][0.8][0]{}
\psfrag{GHL}[tl][tl][0.8][0]{}
$$\includegraphics[width=4cm]{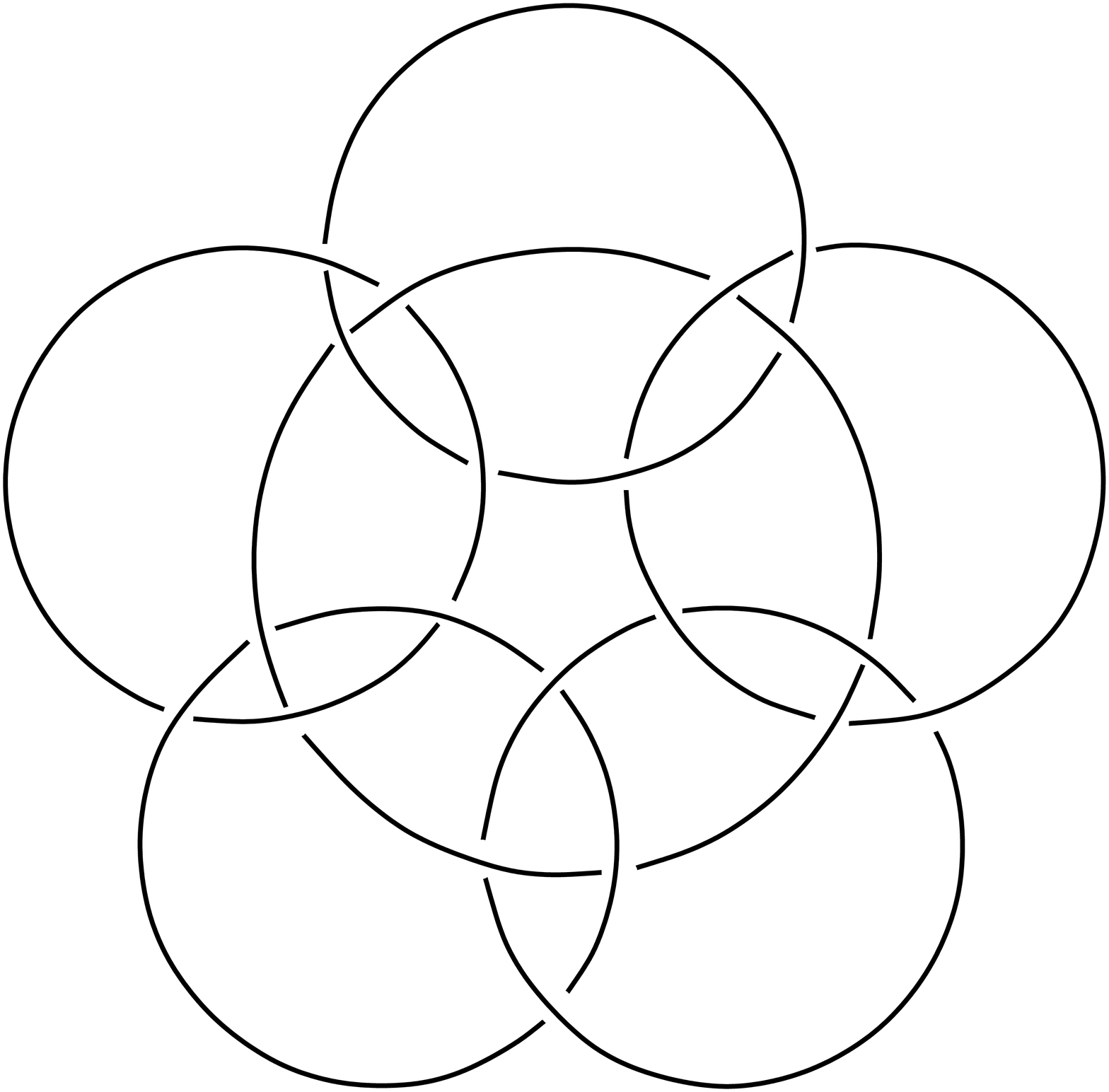}$$
}
Using the techniques of Kawauchi \cite{Kawauchi} one can construct
many other interesting representations $B_L \to \Sigma_n$. 
Unfortunately, Kawauchi's method can not be used to realise representations 
$B_L \to \Sigma_n^+$ that are not conjugate to representations
into $\Sigma_n$ because of some technical restrictions to his 
imitation theory.  Sakuma \cite{Sakuma1} has another interesting family of
examples where $\Zed_{2n} = B_L \to \Sigma_n^+$ is a signed $n$-cycle,
ie: the image is generated by $(123\cdots n-)$, provided $n$ is odd.
{
\psfrag{n}[tl][tl][0.8][0]{}
\psfrag{GHL}[tl][tl][0.8][0]{}
$$\includegraphics[width=4cm]{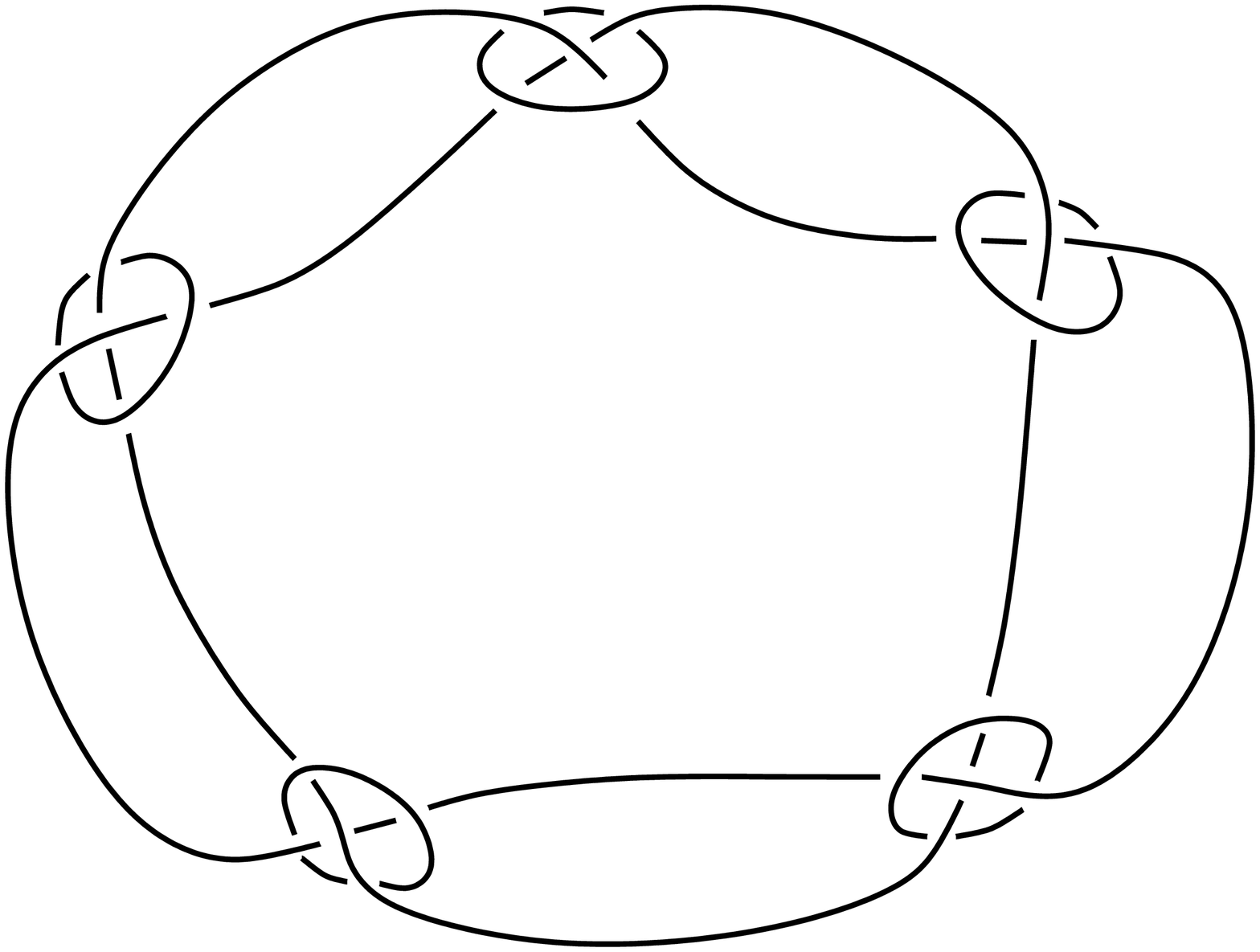}$$
}
There are of course many further directions one might want to
pursue. For a suitable model-space for $\K$, 
$\K$ admits two commuting $SO_2$-actions such that
there is a homotopy-equivalence 
$SO_2\backslash (SO_4 \times \K)/SO_2 \simeq \Emb(S^1,S^3)/\Diff^+(S^1)$.
A natural `next step' would be to determine the $SO_2 \times SO_2$
homotopy type of $\K$. 

\providecommand{\bysame}{\leavevmode\hbox to3em{\hrulefill}\thinspace}

\Addresses

\end{document}